\documentclass[journal]{IEEEtran}

\usepackage{amssymb}
\usepackage[ruled,vlined,linesnumbered]{algorithm2e}
\usepackage{amsthm}
\usepackage{amsfonts}
\usepackage{amsmath}
\usepackage{booktabs}
\usepackage{bm}
\usepackage{color}
\usepackage{cite}
\usepackage{comment}
\usepackage{epsfig}
\usepackage{euscript}
\usepackage{fancybox}
\usepackage{graphicx}
\PassOptionsToPackage{hyphens}{url}\usepackage{hyperref}
\usepackage{lipsum}
\usepackage{pgf}
\usepackage{psfrag}
\usepackage{pifont}
\usepackage{mathrsfs}
\usepackage{multirow}
\usepackage{setspace}
\usepackage{stfloats}
\usepackage{tabularx}
\usepackage{textcomp}
\usepackage[caption=false, font=footnotesize]{subfig}

\newcolumntype{C}[1]{>{\centering\arraybackslash}p{#1}}

\SetKwFunction{Dothat}{Do that}
\SetKw{Kwand}{and}

    \def\Complex{{\rm\rule[.23ex]{.03em}{1.1ex}\kern-.3em{C}}}

    \newcommand{\be}{\begin{equation}} \newcommand{\ee}{\end{equation}}
    \newcommand{\bea}{\begin{eqnarray}} \newcommand{\eea}{\end{eqnarray}}
    \newcommand{\benum}{\begin{enumerate}} \newcommand{\eenum}{\end{enumerate}}

    \newcommand{\qzero}{{\bf 0}}

    \newcommand*{\argmin}{\operatornamewithlimits{argmin}\limits}

\begin{document}

\title{ Placement and Routing Optimization for Automated Inspection with UAVs: A Study in Offshore Wind Farm}

\author{ \IEEEauthorblockN{Hwei-Ming Chung, \IEEEmembership{Student Member, IEEE}, Sabita Maharjan, \IEEEmembership{Senior Member, IEEE}, Yan Zhang, \IEEEmembership{Fellow, IEEE}, Frank Eliassen, \IEEEmembership{Member, IEEE}, and Kai Strunz}

\thanks{This work was supported by Norwegian Research Council under Grants 275106 (LUCS project), 287412 (PACE project), and 267967 (SmartNEM project).
This work was also supported by German Federal Ministry for Economic Affairs and Energy (BMWi) project WindNODE: Echtzeitlabor Energiewende under Grant 03SIN539.
}

\thanks{H.-M. Chung and F. Eliassen are with the Department of Informatics, University of Oslo, Oslo 0373, Norway, (e-mail: hweiminc@ifi.uio.no, frank@ifi.uio.no).}

\thanks{S. Maharjan and Y. Zhang are with Department of Informatics, University of Oslo, Oslo 0373, Norway; and Simula Metropolitan Center for Digital Engineering, Oslo 0167, Norway. (e-mail: yanzhang@ieee.org, sabita@ifi.uio.no).
}

\thanks{K. Strunz is with the Department of Energy and Automation Technology, Technische Universit{\"a}t Berlin, Berlin 10587, Germany (e-mail: kai.strunz@tu-berlin.de).
}

}

\markboth{IEEE Transactions on Industrial Informatics.,~Vol.~XX, No.~X, month ~ year}
{}

\maketitle

\begin{abstract}

Wind power is a clean and widely deployed alternative to reducing our dependence on fossil fuel power generation.
Under this trend, more turbines will be installed in wind farms.
However, the inspection of the turbines in an offshore wind farm is a challenging task because of the harsh environment (e.g., rough sea, strong wind, and so on) that leads to high risk for workers who need to work at considerable height.
Also, inspecting increasing number of turbines requires long man hours.
In this regard, unmanned aerial vehicles (UAVs) can play an important role for automated inspection of the turbines for the operator, thus reducing the inspection time, man hours, and correspondingly the risk for the workers. 
In this case, the optimal number of UAVs enough to inspect all turbines in the wind farm is a crucial parameter.
In addition, finding the optimal path for the UAVs' routes for inspection is also important and is equally challenging.
In this paper, we formulate a placement optimization problem to minimize the number of UAVs in the wind farm and a routing optimization problem to minimize the inspection time.
Wind has an impact on the flying range and the flying speed of UAVs, which is taken into account for both problems. 
The formulated problems are NP-hard.
We therefore design heuristic algorithms to find solutions to both problems, and then analyze the complexity of the proposed algorithms. 
The data of the Walney wind farm are then utilized to evaluate the performance of the proposed algorithms.
Simulation results clearly show that the proposed methods can obtain the optimal routing path for UAVs during the inspection.

\end{abstract}


\section*{Nomenclature} 
\addcontentsline{toc}{section}{Nomenclature}
\begin{IEEEdescription}[\IEEEusemathlabelsep\IEEEsetlabelwidth{$V_1,V_2$}]
\item[A. Sets and Indices]
\item[$i, j$] UAV index.
\item[$k, l$] Turbine index.
\item[$T$]  Total number of turbines in the wind farm.
\item[${\cal N}_{i}$] Set of turbines for UAV $i$ to inspect.
\item[$N$] Number of candidate UAVs in wind farm.
\item[$\mathcal{W}$] Set of wind data.
\\

\item[B. Variables]
\item[$\mathbf{v}_{i, k, l}$]  UAV velocity for UAV $i$ flying from turbine $k$ to turbine $l$.
\item[$||\mathbf{v}_{i, k, l}||_{2}$ ]  Airspeed of UAV $i$ flying from turbine $k$ to turbine $l$.
\item[$\mathbf{s}_{i, k, l}$]  Resultant velocity for UAV $i$ flying from turbine $k$ to turbine $l$.
\item[$||\mathbf{s}_{i, k, l}||_{2}$ ]  Groundspeed of UAV $i$ flying from turbine $k$ to turbine $l$.
\item[$t_{i, k, l}$ ]  Time for UAV $i$ flying from turbine $k$ to turbine $l$.
\item[$\theta_{i, k, l}^{s, w}$] Angle between $\mathbf{s}_{i, k, l}$ and $\mathbf{w}$.
\item[$\theta_{i, k, l}^{s, v}$] Angle between $\mathbf{s}_{i, k, l}$ and $\mathbf{v}_{i, k, l}$ .
\item[$\mathbf{A}$]  Binary vector indicating the UAV state.
\item[$\mathbf{B}$]  Binary matrix indicating the link between UAVs and turbines.
\item[$\mathbf{C}$]  Binary matrix indicating the communication link between UAVs.
\item[$\mathbf{U}_{i}^{m}$] Binary matrix to indicate the $m$-th route of UAV $i$.
\\

\item[C. Parameters]
\item[$\mathbf{q}_{k}$($\mathbf{q}_{i}$)]  Coordinates of turbine $k$ (UAV $i$).
\item[$\mathbf{w}$]  Wind velocity.
\item[$\theta_{w}^{met}$]  Wind direction obtained from wind data.
\item[$\theta_{w}^{pol}$]  Wind direction represented in polar coordinates.
\item[$w_{s}$]  Wind speed obtained from $||\mathbf{w}||_{2}$.
\item[$u_{i}^{max}$]  Maximum flying speed of UAV $i$.
\item[$t_{i}^{max}$]  Maximum flight time of UAV $i$.
\item[$u_{i}^{wind}$]  Maximum resistance to wind of UAV $i$.
\item[$d_{i, j}$($d_{i, k}$)] Distance between UAV $i$ and UAV $j$ (turbine $k$).
\item[$\rho_{i}$] Flying distance of UAV $i$.
\item[$B_{i}^{\mathbf{w}} (\rho_{i})$] Flying range of UAV $i$ under given $\mathbf{w}$ and $\rho_{i}$.
\item[$Z_{i}$] Intersection of flying range under different $\mathbf{w}$.
\item[$p$] Maximum number of turbines for UAVs to inspect.
\item[$d$] Maximum distance for obtaining a communication link between two UAVs.
\item[$M$] Number of routes to perform turbine inspection.
\item[$\mathbf{D}_{i}$] Adjacency matrix for a graph constructed by UAV $i$ and ${\cal N}_{i}$.
\item[$\epsilon_{v}$] Threshold of hourly average wind speed so that hourly wind gust does not exceed $u_{i}^{wind}$.
\item[$path$] Optimal routing path without considering $t_{i}^{max}$.
\\

\item[D. Operators]
\item[$|\cdot|$]   Cardinality of set.
\item[$||\cdot||_{2}$]   Two norm of a vector.
\item[$\qzero$]   Zero vector.
\item[$\mathbf{C}_{ \{ f, g\} }$]   Submatrix that is taken from $\mathbf{C}$ with the row index in set $f$ and the column index in set $g$.
\\

\item[Other notations are defined in the text.]  

\end{IEEEdescription}

\section{Introduction}
\IEEEPARstart{W}{ith} increasing influence of wind power in the energy ecosystem, the capacity of the global wind power is expected to grow by $60\,\%$ over the next $5$ years \cite{wind-power-forecast}.
The overall capacity of all wind turbines installed in the first half of $2019$ in Europe has reached $4.9$ GW \cite{wind-power-europe}.
This is the same amount as the wind power generation capacity installed for the whole year of $2018$.
However, with a large number of wind turbines, inspecting and maintaining the condition of the turbines becomes a challenging task.
Turbines may suffer failures from different components, such as blades, gearbox, yaw system, and so on \cite{2007-turbine-failure}.
The authors in \cite{2014-DTU-report} reported that a blade failure could result in a downtime of more than seven days.

Advanced sensors have been introduced to monitor the operation and the health condition of wind turbines \cite{2015-fibre-lee,2016-lidar-inspect-UAV,2018-blade-thermal-yang,2017-blade-mmwave-wang}.
For instance, in \cite{2015-fibre-lee}, fiber optic sensors were used to monitor the operation of wind turbines.
Lidar sensors were used to detect cracks on the blade in \cite{2016-lidar-inspect-UAV}.
Other technologies used for monitoring and analyzing the conditions of wind turbines are thermal wave radar \cite{2018-blade-thermal-yang} and millimeter wave imaging \cite{2017-blade-mmwave-wang}, respectively.  
Other researchers have developed algorithms for detecting the failure of turbines.
For instance, in \cite{2017-gearbox-failure-DNN}, the authors utilized a deep neural network (DNN) to detect the failure of the gearbox.
The convolutional neural network (CNN) was introduced to detect the icing on the blade in \cite{2019-blade-ice}.

From \cite{2015-fibre-lee,2016-lidar-inspect-UAV,2018-blade-thermal-yang,2017-blade-mmwave-wang,2017-gearbox-failure-DNN,2019-blade-ice}, we can observe that the detection mechanisms are rather advanced.
However, the main challenge is to optimally place the sensors \cite{2015-fibre-lee,2016-lidar-inspect-UAV,2018-blade-thermal-yang,2017-blade-mmwave-wang} or to acquire data \cite{2017-gearbox-failure-DNN,2019-blade-ice}.
For the onshore wind farm, the sensors can be installed near the turbines and then transmit the measurements back to the control center through an aggregator.
However, this method is not suitable for an offshore wind farm because the sensors cannot be easily installed on the sea.
Currently, it is necessary to dispatch qualified service personnel to manually inspect turbines in offshore wind farms.
This procedure may take several days to several weeks, requiring intensive and costly efforts.
Moreover, workers may be subjected to the risks associated with climbing wind turbines and working at height.

Automated inspection of turbines is a solution that can address these issues.
In this case, the unmanned aerial vehicle (UAV) can play a crucial role.
UAVs have been widely applied for the automated inspection for energy systems \cite{2015-UAV-PV-panel,2017-UAV-blade-wang,2017-power-reading-UAV,2019-UAV-blade-wang,2018-UAV-power-inspect,2020-UAV-insulator-inspect,2019-UAV-power-inspect}.
If there is dust on the solar panel, the power generation of the solar panels can be influenced.
The authors in \cite{2015-UAV-PV-panel} proposed a framework for utilizing UAVs to monitor the condition of solar panels.
Regarding wind turbine blades, cracks on the surface can be detected with the help of images taken by the UAVs \cite{2017-UAV-blade-wang,2019-UAV-blade-wang}.
Automatic meter reading is another application of UAVs in power systems, that was studied in \cite{2017-power-reading-UAV}.
The damage of a power line can also be detected by UAVs \cite{2018-UAV-power-inspect} so that workers do not have to climb the transmission line tower.
Moreover, defects of power insulators can be detected by images captured by UAVs and CNN \cite{2020-UAV-insulator-inspect}.
In \cite{2019-UAV-power-inspect}, the authors combined UAVs with fault indicators.
If there is a damage on the distribution line, fault indicators can send a signal to the UAV, and then the UAV can help to transmit the signal back to the operator.
This is similar in wireless communication where the UAV can be regarded as a relay to transmit the signal \cite{2018-UAV-rui-tsg,2019-UAV-rui-twc-secure,2019-UAV-rui-twc-energy}.

Most of the studies mentioned above focused on how to route UAVs to collect the data.
That is, there are several targets for UAVs to collect data from, and therefore the UAV has to find the optimal path to route.
For example, the authors in \cite{2018-UAV-rui-tsg} jointly minimized usage of time and bandwidth to find an optimal path for routing.
Then, maximizing the average secrecy rate to secure the communication link between users and UAVs was considered in \cite{2019-UAV-rui-twc-secure}.
The authors in \cite{2019-UAV-rui-twc-energy} further minimized energy consumption while ensuring that the throughput requirements of the users are met.
However, the placement of UAVs can significantly influence the routing results.
Therefore, the placement problem such as finding the optimal number of UAVs and deriving the optimal topology for UAV placement should be considered as discussed in \cite{2019-UAV-deploy-zhang,2019-UAV-deploy-sun,2018-deplpoy-UAV-predictive}.
The authors in \cite{2019-UAV-deploy-zhang} proposed two algorithms to minimize the time for UAVs flying to the specific locations to serve mobile users.
Specifically, the first algorithm minimized the flight time, and the second algorithm incorporates fairness of allocating transmitting power to users while minimizing the flight time.
In \cite{2019-UAV-deploy-sun}, $K$-means clustering was applied to find the optimal locations and the number of UAVs such that the desired area can be covered with a minimum number of UAVs.
In \cite{2018-deplpoy-UAV-predictive}, the authors applied an algorithm to predict the future throughput requirements of the users.
Then, the optimal number and the topology of UAVs are determined based on the prediction.

In this paper, UAVs are adopted to perform automated monitoring of wind farms.
Specifically, UAVs equipped with sensors, such as Lidar, millimeter wave, or thermographic sensors, can monitor the surface condition of turbines.
Then, we address the placement problem of UAVs in an offshore wind farm and find the optimal routing path for turbine inspection by utilizing UAVs.
Some related research works have focused on visiting targets utilizing multiple UAVs \cite{2018-UAV-rui-tsg,2019-UAV-rui-twc-secure,2019-UAV-rui-twc-energy} and placing UAVs for the optimal topology \cite{2019-UAV-deploy-zhang,2019-UAV-deploy-sun,2018-deplpoy-UAV-predictive}.
However, the placement and routing problems of UAVs in offshore wind farm have not been studied.
Especially, the random realization of an extreme weather event in the offshore wind farm was not considered in those studies.
That is, compared to the onshore situation, wind in offshore wind farms is known to be stronger on average.
The influence of wind speed and wind direction in offshore wind farms is considerable and is important to consider for the UAV placement and routing problems.
The UAVs may crash if the wind speed is over the maximum wind speed resistance of the UAVs during turbine inspection.
Also, wind speed and wind direction impact the flying speed and flying range of UAVs.
We then introduce a mathematical model that addresses the relation between wind speed, wind direction, and UAVs.
This model is essential to consider when formulating the placement and the routing problems for offshore wind farms.
The formulated problems are NP-hard, and therefore we design heuristic algorithms to solve both problems that also take the mathematical model into account.

Overall, the main contributions of this paper are threefold:
\begin{itemize}
\item We present a novel framework for inspecting the wind turbines in the offshore wind farm by utilizing UAVs. 
Then, two optimization problems are formulated for the placement and routing problems, incorporating wind as it can considerably influence the flying range and the flying speed of UAVs.

\item The formulated problems are NP-hard such that they cannot be solved directly.
Therefore, we design heuristic algorithms to obtain the required number and the topology of UAVs in an offshore wind farm and the optimal path for the inspection.

\item We also analyze the complexity of the proposed algorithms.
With the proposed methods, the placement problem can be solved in polynomial time, and the routing problem can be solved with lower complexity compared to the brute-force method.

\end{itemize}

\section{System Model}\label{sec:system_model_problem_formulation}


\subsection{Wind Farm and Wind Model}

The total number of turbines in the offshore wind farm is $T$.
The coordinates of the $k$-th turbine are $\mathbf{q}_{k} = [x_{k}, y_{k}]$.
The wind velocity is denoted by $\mathbf{w} = [ w^{x}, w^{y} ]$.
The projection of the wind velocity on the x-axis and the y-axis are $w^{x}$ and $w^{y}$, respectively.
Quantity $w_{s}$ represents the wind speed, which can be calculated as $w_{s} =||\mathbf{w}||_{2}$.
The wind direction in the polar coordinate system is denoted by $\theta_{w}^{pol}$, which can be calculated from $ \theta_{w}^{pol} = \arctan \frac{w^{y}}{w^{x}}$.
The definition of $\theta_{w}^{pol}$ is different from the wind direction of the meteorological measurements.
Therefore, the wind direction in the meteorological measurements is denoted by $\theta_{w}^{met}$.
In the meteorological measurements, $0$, $\frac{\pi}{2}$, $\pi$, and $\frac{3 \pi}{2}$ are used to represent the north, the east, the south, and the west wind, respectively.
Thus, the phase is represented in a clockwise direction.
In the polar coordinate system however, the phase is represented in a counterclockwise direction.
Therefore, $\theta_{w}^{pol} $ and $\theta_{w}^{met}$ are related as
\begin{equation}
\theta_{w}^{pol} = \frac{3 \pi}{2} - \theta_{w}^{met}.
\end{equation}

\subsection{UAV Model}

In an offshore wind farm, we place UAVs to monitor the condition of the turbines.
Each UAV should inspect the wind turbines assigned to it.
UAV $i$ will be placed at $\mathbf{q}_{i} = [ x_{i}, y_{i} ]$, and the set of the turbines assigned to UAV $i$ is denoted by ${\cal N}_{i}$.
Then, $|{\cal N}_{i}|$ represents the cardinality of set ${\cal N}_{i}$.

\begin{figure}
\begin{center}
\resizebox{2.4in}{!}{%
\includegraphics*{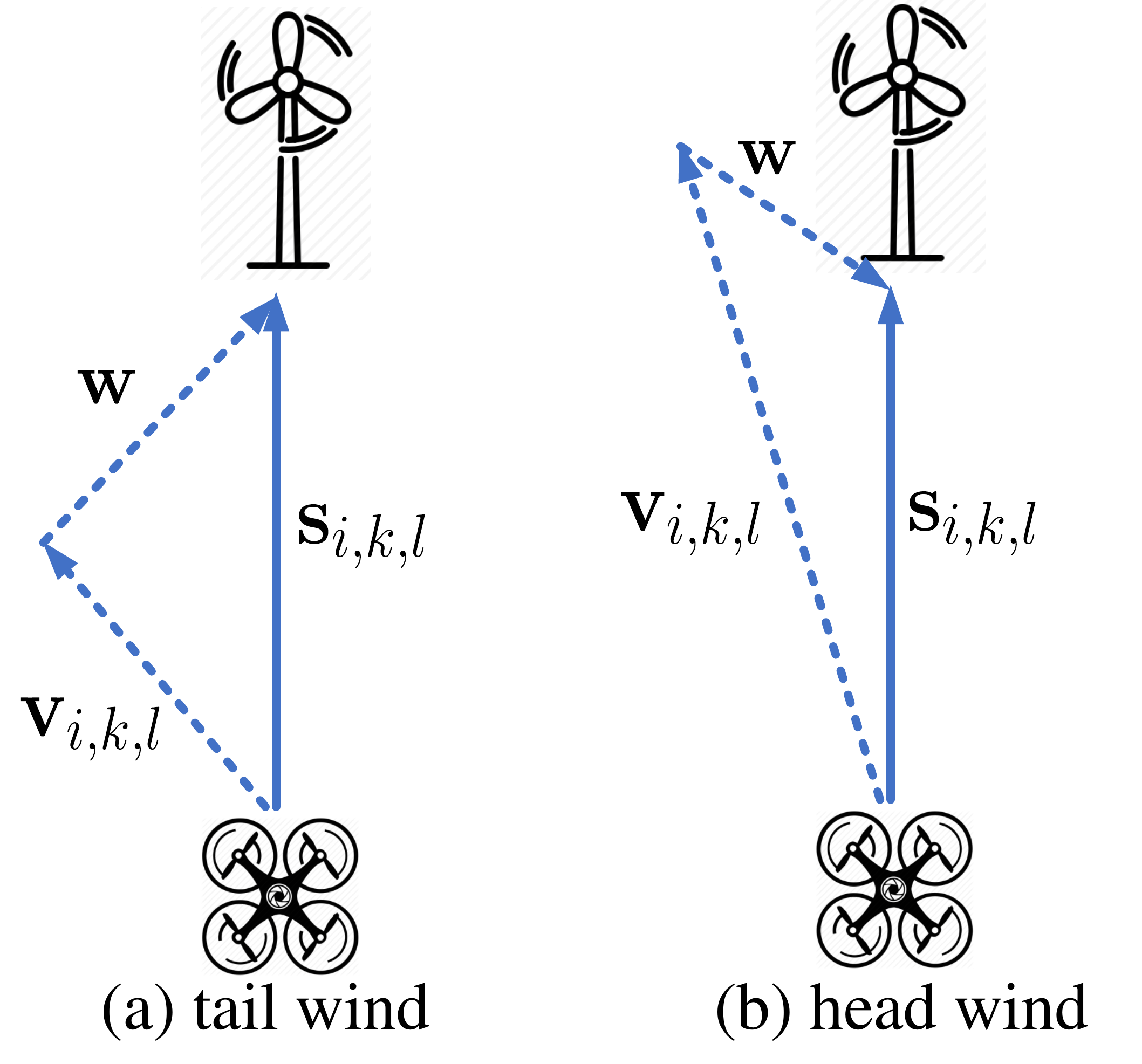} }%
\caption{The relation between UAV and wind.} 
\label{fig:wind_model}
\end{center}
\end{figure}

When UAV $i$ flies to inspect a wind turbine, it may face two wind conditions, namely tail wind and head wind, as shown in Fig. \ref{fig:wind_model}.
The condition of the wind is to be considered in the decision-making of the UAV.
We define $\mathbf{s}_{i, k, l} = [s_{i, k, l}^{x}, s_{i, k, l}^{y} ]$ and $\mathbf{v}_{i, k, l} = [ v_{i, k, l}^{x}, v_{i, k, l}^{y} ]$ as the resultant velocity and the UAV velocity of UAV $i$ flying from turbine $k$ to turbine $l$, respectively.
Components $s_{i, k, l}^{x}$ and $v_{i, k, l}^{x}$ are the projections on the x-axis, and $s_{i, k, l}^{y}$ and $v_{i, k, l}^{y}$ are the projections on the y-axis.
The UAV velocity is the initial velocity of the UAV, and the resultant velocity is the velocity influenced by the wind.
The relation between the UAV velocity, the wind, and the resultant velocity is expressed as
\begin{equation} \label{eq:velocity-relation}
\mathbf{v}_{i, k, l} + \mathbf{w} = \mathbf{s}_{i, k, l}.
\end{equation}
Quantities $||\mathbf{v}_{i, k, l}||_{2}$ and $||\mathbf{s}_{i, k, l}||_{2}$ are regarded as airspeed and groundspeed, respectively.
UAV $i$ has the maximum speed limit of $u_{i}^{max}$.
Usually,  $u_{i}^{max}$ is referred to as the maximum value for the airspeed.
However, UAVs may not remain stable, and the structural capacity of UAVs may degrade if UAVs fly at too high groundspeed.
Therefore, airspeed and groundspeed are both limited to $u_{i}^{max}$ in this paper.
Specifically, for the tail wind situation in Fig. \ref{fig:wind_model}(a), the groundspeed is limited to this value.
Then, when the UAV faces a head wind, the airspeed is limited to $u_{i}^{max}$.
The angle between $\mathbf{s}_{i, k, l}$ and $\mathbf{v}_{i, k, l}$ is denoted by $\theta_{i, k, l}^{s, v}$.
Also, $\theta_{i, k, l}^{s, w}$ is used to represent the angle between $\mathbf{w}$ and $\mathbf{s}_{i, k, l}$.
For UAV $i$, the maximum wind speed resistance is denoted by $u_{i}^{wind}$.

The time for UAV $i$ to travel from turbine $k$ to turbine $l$ can be calculated as
\begin{equation} \label{eq:travel_time}
t_{i, k, l} = \frac{ ||\mathbf{q}_{l} - \mathbf{q}_{k}||_{2} }{||\mathbf{s}_{i, k, l}||_{2}}.
\end{equation}
There exists also a maximum flight time for UAV $i$ denoted by $t_{i}^{max}$, which represents an upper limit of the total flight time during the inspection.
The distance between UAV $i$ and UAV $j$ is denoted by $d_{i, j}$, and the distance between UAV $i$ and turbine $k$ is represented as $d_{i, k}$.

\section{Problem Formulation}\label{sec:problem_formulation}

In this section, we formulate a placement optimization problem to obtain the required number of UAVs in the wind farm and the corresponding topology.
Then, a routing optimization problem is formulated to find the optimal path for the inspection.
Wind is a very important factor in determining the flying range and flying speed of UAVs.
Therefore, the influence of the wind is incorporated in both problems.

\subsection{Placement Problem Formulation} \label{subsec:deploy}

We introduce an important parameter of UAV $i$, namely the flying range.
Here, the UAVs are assumed to fly at the same altitude, and therefore the z-axis can be ignored.
The flying range of UAV $i$ under the wind condition $\mathbf{w}$ can be expressed as 
\begin{equation} \label{eq:flying_range}
B_{i}^{\mathbf{w}} (\rho_{i}) =  \{ x, y \in \mathbb{R}: ||\mathbf{r}||_{2} \leq \rho_{i} \},
\end{equation}
where
\begin{equation}
\mathbf{r} = [ x-x_{r}, y - y_{r} ].
\end{equation}
Quantity $\rho_{i}$ is the actual flying distance of UAV $i$, which can be calculated by
\begin{equation} \label{eq:flying_range_max}
\rho_{i} = \frac{u_{i}^{max}t_{i}^{max}}{2}, 
\end{equation}
assuming that the UAV can fly from and then back to the starting point at maximum flying speed, $u_{i}^{max}$, during maximum flight time, $t_{i}^{max}$.
However, $t_{i}^{max}$ is undetermined to the UAV because it is influenced by $u_{i}^{max}$.
To obtain the value of $\rho_{i}$, an analytical model is introduced in Appendix \ref{subsec:UAV_energy_model} to determine $u_{i}^{max}$ and $t_{i}^{max}$.
The flying range of the UAV is regarded as a circle, with $x_{r}$ and $y_{r}$ as the center of the circle, which is calculated as 
\begin{equation} \label{eq:center_point_cal}
\left\{
\begin{array}{ll}
                x_{r} = x_{i} + w^{x} t_{i}^{max},    \\
                y_{r} = y_{i} + w^{y} t_{i}^{max}.
\end{array}
        \right.
\end{equation}
The flying range of UAV $i$ is influenced by different wind conditions.
Thus, the flying range of the UAV is the intersection of the flying range under different wind conditions, i.e., 
\begin{equation}\label{eq:flying_range_intersec}
Z_{i} = \bigcap_{\mathbf{w}} B_{i}^{\mathbf{w}} (\rho_{i}).
\end{equation}

Initially, a total of $N$ candidate UAVs are placed in the wind farm.
We then introduce matrices $\mathbf{A} = [a_{i}]_{1\times N}$, $\mathbf{B} = [b_{i, k}]_{N \times T}$, and $\mathbf{C} = [c_{i, j}]_{N \times N}$ to denote the states of candidate UAVs, the association between UAVs and turbines, and the communication link between UAVs, respectively.
Here, $a_{i}$ is set to $1$ when this candidate UAV should be placed in the wind farm, and the candidate UAV should be removed if $a_{i}$ is $0$.
Also, a docking station \cite{2016-UAV-dockstation}\footnote{https://www.airoboticsdrones.com/} is placed together with the UAV so that the UAV can charge or swap its battery.
If turbine $k$ is assigned to UAV $i$, $b_{i, k}$ is set to $1$; otherwise, $b_{i, k}$ is $0$.
When determining the topology of the UAVs, we need to ensure that the UAVs can maintain communication links with each other.
By doing so, the collisions between UAVs can be prevented \cite{2007-UAV-Adhoc-collision,2018-UAV-Adhoc-collision}.
Quantity $c_{i, j}$ is set to $1$ if UAV $i$ and $j$ obtain a communication link.

The objective function is to minimize the number of UAVs that should be placed in the wind farm.
Therefore, the placement optimization problem takes the form:
\begingroup
\allowdisplaybreaks
\begin{subequations}\label{eq:UAV_deploy_problem}
\begin{align}
\min_{ \mathbf{A}, \mathbf{B}, \mathbf{C}, x_{i}, y_{i} } &  ~  \sum_{i=1}^{N} a_{i}  & \\
\mbox{subject to}    & ~~ a_{i}, b_{i, k}, c_{i, j} \in \{ 0, 1\}, & \!\!\!\!\!\forall i, k, j \label{eq:bina_para}\\
					 & ~~ \sum_{i=1}^{N} b_{i, k} \leq 1, & \!\!\!\!\!\forall k  \label{eq:UAV_only1_tur}\\
                     & ~~ \sum_{k=1}^{T} b_{i, k} \leq p, & \!\!\!\!\!\forall i  \label{eq:UAV_max_tur}\\
                     & ~~ [x_{k}, y_{k}] \in Z_{i}  & \!\!\!\!\!\forall b_{i, k} =1 \label{eq:range_turbine}\\
                     & ~~ [x_{i}, y_{i}] \in \{ [x_{k}, y_{k}] \}  & \!\!\!\!\!\forall b_{i, k} =1 \label{eq:UAV_pos}  \\
                     & ~~ b_{i, k} \leq a_{i}, & \forall k \label{eq:only_active_UAV} \\
                     & ~~ \sum_{i=1}^{N} \sum_{k=1}^{M} b_{i, k} \geq T, &  \label{eq:all_turbine_assign}\\
                     & ~~ \sum_{j=1, j \neq i}^{N} c_{i, j} \geq 1, & \!\!\!\!\!\forall a_{i}=1 \label{eq:commu_num}\\
                     & ~~ \sqrt{ (x_{i} - x_{j})^{2} +  (y_{i} - y_{j})^{2} } \leq d, & \!\!\!\!\! \forall i, j. \label{eq:distance_limit}
\end{align}
\end{subequations}
\endgroup
Quantities $a_{i}$, $b_{i, k}$, and $c_{i, j}$ are defined as binary variables in (\ref{eq:bina_para}).
Constraints (\ref{eq:UAV_only1_tur}) and (\ref{eq:UAV_max_tur}) state that each turbine can only be assigned to one UAV, and each UAV can inspect up to $p$ turbines, respectively.
Here, $p$ is defined by the wind farm operator.
The location $[x_{k}, y_{k}]$ of turbine $k$ assigned to UAV $i$ must be in the flying range of UAV $i$ as given by (\ref{eq:range_turbine}).
The UAVs cannot be placed on the sea, and therefore each UAV should be placed inside one of the turbines assigned to it as indicated in (\ref{eq:UAV_pos}).
Also, constraint (\ref{eq:only_active_UAV}) states that the turbine can only be assigned to a UAV that is actually placed in the wind farm.
Constraint (\ref{eq:all_turbine_assign}) indicates that all $T$ turbines must be assigned to the active UAVs.
The minimum number of communication links one UAV must have is given by (\ref{eq:commu_num}).
Finally, (\ref{eq:distance_limit}) indicates that the distance between any two UAVs should be lower than $d$.

\subsection{Routing Problem Formulation} \label{subsec:route}

For a given topology, we study how to route the UAVs to inspect the wind turbines.
A UAV and the turbines assigned to it can be represented in a graph, defined as $\mathcal{G}_{i}=\{ {\cal N}_{i}, {\cal E}_{i}$\}, where ${\cal N}_{i}$ is the set of turbines assigned to UAV $i$, and will be the nodes in the graph, and then ${\cal E}_{i}$ is the set of the edges which connects each turbine.
With the graph structure, we can create an adjacency matrix for the graph denoted by $\mathbf{D}_{i}$.
The value of the $k$-th column and the $l$-th row in $\mathbf{D}_{i}$ is $t_{i, k, l}$ as the flight time from turbine $k$ to $l$.
The value of $t_{i, k, l}$ and the value of $t_{i, l, k}$ are not the same because of the wind, and therefore the adjacency matrix is not symmetric.
Thus, $\mathcal{G}_{i}$ is an asymmetric graph.

In the routing problem, $M$ denotes the number of required routes to inspect the turbines.
We introduce another matrix $\mathbf{U}_{i}^{m} = [U_{i, k, l}^{m}]_{|{\cal N}_{i}|\times |{\cal N}_{i}|}$ to denote the $m$-th route for UAV $i$.
More specifically, $U_{i, k, l}^{m}$ is $1$ when the UAV chooses to fly from turbine $k$ to turbine $l$; otherwise, $U_{i, k, l}^{m}$ is $0$.
The routing problem can then be formulated as
\begingroup
\allowdisplaybreaks
\begin{subequations}\label{eq:UAV_routing_problem}
\begin{align}
\min_{\substack{ M, \mathbf{U}_{i}^{m}, \mathbf{v}_{i, k, l},\\ \mathbf{s}_{i, k, l}, \theta_{i, k, l}^{s, v}}}  &  ~ \sum_{m=1}^{M}\sum_{k \in {\cal N}_{i} } \sum_{l \in {\cal N}_{i} \setminus \{ k\} } t_{i, k, l} U_{i, k, l}^{m}   \\
\mbox{subject to}    & ~~~~ U_{i, k, l}^{m} \in \{ 0, 1\}, \quad\qquad\qquad\forall k, l \in  {\cal N}_{i} \label{eq:binary_limit} \\
					 & ~~ \sum_{k \in {\cal N}_{i}} U_{i, s, k}^{m} = \sum_{k \in {\cal N}_{i}} U_{i, k, s}^{m} = 1, \quad\quad \forall m \label{eq:start_point} \\
                     & ~~ \sum_{l \in {\cal N}_{i} \setminus \{ l\}  } U_{i, l, k}^{m} \!=\! \sum_{l \in {\cal N}_{i} \setminus \{ k\}  } U_{i, k, l}^{m} = 1, \forall m \label{eq:link_point} \\
                     & ~~ \sum_{k \in Q} \sum_{l \in Q} U_{i, k, l}^{m} \leq |Q|\!-\! 1,   \!\forall Q  \subsetneq  {\cal N}_{i}, m, |Q|>2 \label{eq:avoid_subtour} \\
                     & ~~ 1 \leq M \leq |{\cal N}_{i}|-1 \label{eq:route_limit}  \\
                     & ~~ \sum_{k, l \in {\cal N}_{i}} t_{i, k, l} U_{i, k, l}^{m} \leq t_{i}^{max}, ~\quad\qquad\quad\forall m \label{eq:time_limit}\\
                     & ~~ ||\mathbf{v}_{i, k, l}||_{2}\leq u_{i}^{max},  \,\forall k, l \in  {\cal N}_{i} , U_{i, k, l}^{m} = 1 \label{eq:airspeed_limit}\\
                     & ~~ ||\mathbf{s}_{i, k, l}||_{2}\leq u_{i}^{max},  ~\forall k, l \in  {\cal N}_{i} , U_{i, k, l}^{m} = 1 \label{eq:ground_limit}\\
                     & ~~ \mathbf{v}_{i, k, l} + \mathbf{w} = \mathbf{s}_{i, k, l},   \forall k, l \in  {\cal N}_{i} , U_{i, k, l}^{m} = 1 \label{eq:wind_speed_relation}
\end{align}
\end{subequations}
\endgroup
In (\ref{eq:UAV_routing_problem}), the objective is to minimize the flight time and the number of routes for inspecting the turbines.
Here, $U_{i, k, l}^{m}$ is a binary parameter as shown in (\ref{eq:binary_limit}) to represent the path of routing during the inspection.
Eq. (\ref{eq:start_point}) indicates that the starting point of every route should be $s$, which is the position of the UAV, $\mathbf{q}_{i}$.
Then, there can only exist one route between turbines as stated in (\ref{eq:link_point}).
Constraint (\ref{eq:avoid_subtour}) ensures that a closed path does not exist in the subset $Q$ of ${\cal N}_{i}$.
Constraint (\ref{eq:route_limit}) states that the number of routes must be less than the number of turbines in ${\cal N}_{i}$.
The summation of the flight times in every route is not to exceed $t_{i}^{max}$ according to (\ref{eq:time_limit}).
Constraints (\ref{eq:airspeed_limit}) and (\ref{eq:ground_limit}) enforce that the airspeed and groundspeed are bounded by the maximum speed, respectively.
The relationship between the wind, the UAV velocity, and the resultant velocity is given by (\ref{eq:wind_speed_relation}).

\section{Algorithm Design} \label{sec:algorithm_design}

The formulations in (\ref{eq:UAV_deploy_problem}) and (\ref{eq:UAV_routing_problem}) cannot be solved directly as they contain binary parameters.
Both problems are mixed-integer linear programming (MILP) problems.
In this case, if the dimension of the problems increases, the problems may become NP-hard.
We therefore design heuristic algorithms to solve the problems.

\subsection{Flying Range Determination}

For both placement and routing optimization problems, an important parameter is the flying range.
The flying range can be influenced by the wind conditions in the wind farm.
However, wind conditions in the future are unknown to UAVs, and different UAVs have different maximum wind speed resistance.
Since the algorithm is used to determine how many UAVs are required for inspections and where to place the UAVs, the algorithm can be considered as part of a planning stage before operation begins.
As such, the algorithm makes use of historic wind data.

The wind data, denoted by $\mathcal{W}$, provide the hourly average wind speeds and directions for several days.
Selected data also provide the hourly wind gust and the corresponding direction.
The UAV cannot perform inspection if the wind speed exceeds maximum wind speed resistance, $u_{i}^{wind}$, in that hour.
Therefore, we introduce an auxiliary parameter, $\epsilon_{v}$.
This parameter is a user input.
Specifically, an appropriate value is assigned to $\epsilon_{v}$ such that, based on a stochastic analysis of the historic data, the hourly wind gust does not exceed $u_{i}^{wind}$ if the hourly average wind speed is lower than $\epsilon_{v}$.
The assignment of $\epsilon_{v}$ is done prior to the start of Algorithm \ref{ago:range_deter}.
For usage of $\epsilon_{v}$ in the loop starting in line $2$, the angle between $0$ and $2 \pi$ is discretized into $\mu$ segments of equal size.
For every segment, the wind velocity located in this segment is taken out to construct a subset $\mathcal{W}_{b}$.
Then, we compare the maximum hourly average wind speed in $\mathcal{W}_{b}$ with $\epsilon_{v}$.
The smallest value is retained to define the flying range according to (\ref{eq:flying_range}).
After executing the loop for all segments, the flying range for the UAV can be obtained from (\ref{eq:flying_range_intersec}).
The detailed steps are summarized in Algorithm \ref{ago:range_deter}.

\begin{algorithm}
\small
\caption{Flying Range Determination}
\label{ago:range_deter}
 \DontPrintSemicolon
\KwIn{wind data $\mathcal{W}$, $\epsilon_{v}$, $t_{i}^{max}$, $u_{i}^{max}$}
\KwOut{$Z_{i}$}
Calculate $\rho_{i}$ based on $t_{i}^{max}$ and $u_{i}^{max}$ using (\ref{eq:flying_range_max})\;
\For{$b  = 1$ \KwTo $\mu$ }{
	$\mathcal{W}_{b} = \left \{ \mathbf{w}  \mid \mathbf{w} \in \mathcal{W}, \frac{2\pi \times (b-1)}{\mu} \leq \theta_{w}^{pol} \leq \frac{2\pi \times b}{\mu} \right\}$ \;
	$ \hat{w_{s}} = \min \{ \epsilon_{v}, \max_{ \mathbf{w} \in \mathcal{W}_{b}} ||\mathbf{w}||_{2}\}$ \;
	$ \hat{\theta_{w}^{pol}} = \frac{2b\pi-\pi}{\mu}$ \;
	$ \hat{\mathbf{w}} = \left[ \hat{w_{s}} \cos \left(\hat{\theta_{w}^{pol}} \right), \hat{w_{s}} \sin \left( \hat{\theta_{w}^{pol}} \right) \right]$ \;
	Calculate the flying range with (\ref{eq:flying_range}) and $\hat{\mathbf{w}}$\;
   }
Perform the intersection as mentioned in (\ref{eq:flying_range_intersec})\;
\end{algorithm}

\subsection{Algorithms for Obtaining Topology of UAVs}

With known flying range, we can design an algorithm to get the topology of UAVs in the wind farm.
At the beginning of the search, $T$ UAVs are placed in the wind farm; in this case, every turbine has a UAV assigned to it.
If the distance between two UAVs is shorter than $d$, they establish a communication link.
Then, the turbines inside the flying range of UAV $i$ are assigned to UAV $i$.
The detailed process of the initialization is provided in Algorithm \ref{ago:deploy_ini}.

\begin{algorithm}
\small
\caption{Initialization}
\label{ago:deploy_ini}
 \DontPrintSemicolon
\KwIn{The coordinates of $T$ turbines, $Z_{i}$}
\KwOut{$\mathbf{A}$, $\mathbf{B}$, $\mathbf{C}$, $x_{i}$, $y_{i}$}
Set $N = T$, set every $[x_{i}, y_{i}]$ to $ [x_{k}, y_{k}]$ \;
\For{$i = 1$ \KwTo $N$ \Kwand $j = 1$ \KwTo $N$}{
    $a_{i} = 1$ \;
    Calculate the distance between UAV $i$ and $j$, $d_{i, j}$\;
    \eIf {$d_{i, j} \leq d$}{
    $c_{i, j} = 1$ \;
    }
    {
    $c_{i, j} = 0$ \;
    }
}
\For{$i = 1$ \KwTo $N$ \Kwand $k = 1$ \KwTo $T$}{
    \eIf {$[x_{k}, y_{k}] \in Z_{i}$}{
    $b_{i, k} = 1$ \;
    }
    {
    $b_{i, k} = 0$ \;
    }
}
\end{algorithm}

After initialization, constraint (\ref{eq:UAV_max_tur}) should be validated.
That is, some UAVs may have to inspect more than $p$ turbines.
Redundant connections between the UAV and the turbines assigned to it should be deleted.
For that, the distance between the UAV and the turbines assigned to it is sorted in a decreasing order.
Then, turbines in ${\cal N}_{i}$ are reassigned if UAV $i$ is assigned more than $p$ turbines; otherwise, we go to next UAV.
The deletion starts from the turbine which has the longest distance to UAV $i$.
At the same time, we have to ensure that this turbine is inspected by another UAV. 
Then, the procedure of deleting the turbines is repeated until $|{\cal N}_{i}| \leq p$.
The details are provided in Algorithm \ref{ago:cons_check}.

\begin{algorithm}
\small
\caption{Restrict Inspection Limit of UAVs}
\label{ago:cons_check}
 \DontPrintSemicolon
\KwIn{$\mathbf{B}$, $p$}
\KwOut{$\mathbf{B}$}
\For{$i = 1$ \KwTo $N$}{
Sort turbines in ${\cal N}_i$ based on their distance to UAV $i$ in a decreasing order as $e_1,e_2,\dots ,e_{|{\cal N}_i|}$\;
$k = e_{1}$\;
\While{$|{\cal N}_{i}| > p$}{
    \eIf{$\sum_{j, j \neq i} b_{j, k} \geq 1$}{
        Set $b_{i, k} = 0$ and remove $k$ from ${\cal N}_{i}$\;
        }
        {
        Set $k$ to next turbine \; 
        }
}
}
\end{algorithm}

After executing Algorithm \ref{ago:cons_check}, the current solution satisfies constraints (\ref{eq:bina_para}) and (\ref{eq:UAV_max_tur})-(\ref{eq:distance_limit}).
However, some turbines may be assigned to more than one UAV.
Also, the current number of placed UAVs is not minimized and still equal to the number of turbines, i.e., $\sum_{i} a_{i} = T$.
Therefore, we need to reduce the number of the placed UAVs and then fix the issue of one turbine being assigned to multiple UAVs given by (\ref{eq:UAV_only1_tur}).
The details of this process are summarized in Algorithm \ref{ago:UAV_num_mininize}.
In this algorithm, two auxiliary parameters are introduced, namely $cur$ and $step$.
Parameter $cur$ is a set to represent the set of the UAVs with $a_{i}=1$.
Then, $step$ indicates the current iteration step, and it will be used in the simulation part.
In line $2$, the UAVs in $cur$ are sorted based on the number of turbines intersecting with other UAVs in a decreasing order as $f_1,f_2,\dots ,f_{|{cur}|}$.
Then, $i$ is set to $f_{1}$ in line $4$.
In lines $5$ and $6$, we need to find if any turbine within ${\cal N}_{i}$ can also be served by other UAVs.
If that is true, UAV $i$ is removed as described in lines $7$ to $9$.
Otherwise, the connection between the UAV and the turbine is deleted based on the distance in lines $10$ to $16$.

\begin{algorithm}
\small
\caption{UAV Number Minimization}
\label{ago:UAV_num_mininize}
\DontPrintSemicolon
\KwIn{$\mathbf{A}$, $\mathbf{B}$, $\mathbf{C}$}
\KwOut{$\mathbf{A}$, $\mathbf{B}$, $\mathbf{C}$}
$step = 0$, $cur = \{ 1, 2, \dots, N\}$ \;
Sort UAV $i$ in $cur$ based on $\sum_{j\in cur\setminus\{i \} }|{\cal N}_{i}  \cap {\cal N}_{j}|$ in a decreasing order as $f_1,f_2,\dots ,f_{|{cur}|}$\;  
\While{$|{\cal N}_{i}  \cap {\cal N}_{j}|>0$}{
    $diff = \{ \}$, $delete =$ True,  $i = f_{1}$\;
    \For{$j = 1$ \KwTo $N$ \Kwand $a_{j} = 1$}{
    $diff = diff \cup ({\cal N}_{i}  \cap {\cal N}_{j})$ \;
    }  
    \If{$diff = {\cal N}_{i}$ \!\Kwand \!$\mathbf{C}_{\{j \in cur \setminus \{ i\}, cur \setminus \{ i\} \}} \!\!\neq\!\! \qzero $ \Kwand $delete$}{
            $a_{i} \!=\! b_{i, k} \!=\! 0 ~\forall k \!\in\! {\cal N}_{i}$, $cur \!=\! cur \!\setminus \!\{ i\}\!$, $step \!=\! step \!+\! 1$\;
            $c_{i, j} = c_{j, i} = 0 ~\forall j \in cur $, $delete =$ False \;
        }
    \If{$delete$}{

        \For{$j = f_{1}$ \KwTo $f_{|{cur}|}$}{
            \eIf {$d_{j, k} < d_{i, k}$}{
            $b_{i, k} = 0$ \;
            }
            {
            $b_{j, k} = 0$ \;
            }
        }
        $step = step + 1$ \;
        
    }
    Repeat $2$ \;
    
}
\end{algorithm}

\subsection{Algorithms for Finding Optimal Routing Path}

With the topology of the UAVs, we can now introduce how to route the UAVs to inspect the wind turbines.
In (\ref{eq:UAV_routing_problem}), the optimal values of several parameters should be found.
Moreover, the variables at the upper bounds of summation and the binary parameters make the problem difficult to solve.
To address the challenge, solving (\ref{eq:UAV_routing_problem}) is separated into three stages and each of them is solved individually.

In the first stage, the adjacency matrix, $\mathbf{D}_{i}$, should be constructed.
Prior to calculating $\mathbf{D}_{i}$, we have to calculate $\mathbf{s}_{i, k, l}$, $\mathbf{v}_{i, k, l}$, and $t_{i, k, l}$ for all $k, l \in {\cal N}_{i}$.
The calculation of $\mathbf{s}_{i, k, l}$ and $\mathbf{v}_{i, k, l}$ differs depending on whether the UAV is facing head wind or tail wind.
Quantity $\theta_{i, k, l}^{s, w}$ is utilized to determine the wind condition, and it can be obtained by calculating the inner product of $\mathbf{s}_{i, k, l}$ and $\mathbf{w}$.
If $\theta_{i, k, l}^{s, w}$ is between $0$ and $\frac{\pi}{2}$, the UAV is facing a tail wind; otherwise, the UAV $i$ is facing a head wind.
The resultant velocity of UAV $i$ flying from turbine $k$ to turbine $l$ can be obtained from
\begin{equation} \label{eq:s_veloci_cal}
\mathbf{s}_{i, k, l} \!\!=\!\! \left\{
\begin{array}{lll}
                \left[u_{i}^{max} \cos(\theta_{s}), u_{i}^{max} \sin(\theta_{s}) \right], &  0 \leq \theta_{i, k, l}^{s, w} \leq \frac{\pi}{2}, \\
                \left[ u_{i}^{s} \cos(\theta_{s}), u_{i}^{s} \sin(\theta_{s}) \right], & \frac{\pi}{2} < \theta_{i, k, l}^{s, w} \leq \pi,
\end{array}
        \right.
\end{equation}
where $\theta_{s}$ is given by $\arctan((y_{l}-y_{k})/(x_{l}-x_{k}))$.
In (\ref{eq:s_veloci_cal}), $u_{i}^{s}$ can be calculated using
\begin{equation}
u_{i}^{s} = u_{i}^{max} \cos (\theta_{i, k, l}^{s, v}) - w_{s} \cos (\pi - \theta_{i, k, l}^{s, w}),
\end{equation}
where
\begin{equation} \label{eq:sv_theta}
\theta_{i, k, l}^{s, v} = \arcsin \frac{ w_{s} \sin (\pi - \theta_{i, k, l}^{s, w}) }{ u_{i}^{max} }.
\end{equation}
Then, the relation in (\ref{eq:velocity-relation}) can be used to calculate the UAV velocity of UAV $i$, $\mathbf{v}_{i, k, l}$.
With $\mathbf{s}_{i, k, l}$, $\mathbf{q}_{l}$, and $\mathbf{q}_{k}$, $t_{i, k, l}$ can be obtained by (\ref{eq:travel_time}), and then used to construct $\mathbf{D}_{i}$.

In the next stage, we relax the constraints (\ref{eq:route_limit})-(\ref{eq:wind_speed_relation}) and solve the following optimization problem 
\begin{subequations}\label{eq:routing_problem-1}
\begin{align}
\min  &  ~ \sum_{m=1}^{1}\sum_{k \in {\cal N}_{i} } \sum_{l \in {\cal N}_{i} \setminus \{ k\} } t_{i, k, l}  U_{i, k, l}^{m}   \\
\mbox{subject to}    & ~~\mbox{(\ref{eq:start_point})}-\mbox{(\ref{eq:avoid_subtour})}.
\end{align}
\end{subequations}
In the above problem, the optimal routing paths for the UAVs are searched without considering the time limit.
A heuristic algorithm is designed to solve (\ref{eq:routing_problem-1}), and it is presented as Algorithm \ref{ago:optimal-route}.
At the beginning of the algorithm, the number of routes is set to $1$, and the starting point, $s$, is set to the location of UAV $i$, $[x_{i}, y_{i}]$.
Then, we create three sets, namely $From$, $To$, and $To_{next}$.
The possible set of the current locations of the UAV is defined as $From$.
The sets $To$ and $To_{next}$ denote the possible sets of remaining turbines for a UAV to inspect in the next step and next two steps, respectively.
Sets $From$, $To$, and $To_{next}$ are initialized by ${\cal N}_{i}\setminus\{s\}$, $\{s\}$, and an empty set, respectively.
Another vector, $path$, is used to denote the optimal routing path for the UAV. 
Then, $g (k, To)$ denotes the minimal time to finish the inspection of remaining turbines in the set $To$ when the UAV is at the $k$-th turbine.
At Lines $2$ to $9$ in Algorithm \ref{ago:optimal-route}, the elements in the set $From$ is moved to set $To$, and then we calculate the $g (k, To)$ with $k \in From$.
If $To = \{s\}$, $To_{next} = \{\}$, and the UAV is at turbine $k$, $g (k, To)$ is obtained as
\begin{equation} \label{eq:path_cost_first}
g (k, To) =  t_{i, k, s},
\end{equation}
where the minimal time for finishing the inspection is to fly back to the starting point from turbines $k$.
\begin{equation}\label{eq:path_cost}
g (k, To) = \min_{l \in To, l \neq s, To_{next} = To \setminus \{ l\}}  t_{i, k, l} + g (l, To_{next}).
\end{equation}
With all $g (k, To)$, a path with minimal routing time, $path$, can be obtained as shown in lines $10$ to $15$ in Algorithm \ref{ago:optimal-route}.

\begin{algorithm}
\small
\caption{Search Optimal Routing Path}
\label{ago:optimal-route}
 \DontPrintSemicolon
\KwIn{$\mathbf{D}_{i}$, the coordinate of UAV $i$ $[x_{i}, y_{i}]$}
\KwOut{$path$}
Set $M=1$; $s= [x_{i}, y_{i}]$; $From = {\cal N}_{i}\setminus\{s\}$; $To = \{ s \}$; $To_{next}=\{\}$; $path = [0]_{(|{\cal N}_{i}|+1) \times 1} $\;
\For{$num  = 0$ \KwTo $|{\cal N}_{i}|-2$}{
\For{ $K = 1$ \KwTo $\binom{|{\cal N}_{i}|-1}{num}$ }{
    Move different $num$ elements in the set $From$ to the set $To$\;
    \For{$k \in From$}{
        \eIf {$num =0$}{
            Calculate (\ref{eq:path_cost_first})\;
        }
        {
	        Calculate (\ref{eq:path_cost})\;
        }
    }
}
}
$From = \{ s \}$ ; $To = {\cal N}_{i} \setminus \{ s \}$ ; $path[1]=s$\;
\For{$num  = 2$ \KwTo $|{\cal N}_{i}|$}{
    $k = path[num-1]$ \;
    $l = \argmin_{l \in To, To_{next} = To \setminus \{ l\}}  t_{i, k, l} + g (l, To_{next})$\;
    $path[num]=l$, $From = From \cup \{l\}$, $To = To \setminus \{ l \}$\; 
}
$path [|{\cal N}_{i}|+1] = s $
\end{algorithm}

The results of Algorithm \ref{ago:optimal-route} may yield a total flight time that exceeds the time limit.
In this case, the routing path is modified based on the maximum flight time.
More specifically, the UAV should be able to fly back to the starting point to charge or swap its battery that is low or empty.
Therefore, it is tested if the UAV is able to fly back to $s$ when it decides to inspect the $l$-th turbine from the $k$-th turbine.
In this case, two auxiliary parameters, $t_{accu}$ and $t_{compare}$, are introduced.
Quantity $t_{accu}$ is used to denote the accumulative flight time of flying from $s$ to turbine $l$ via turbine $k$.
Then, the time for the UAV to fly from turbine $l$ back to $s$ is denoted by $t_{compare}$. 
If $t_{accu}$ and $t_{compare}$ are all below $t_{i}^{max}$, UAV $i$ can fly from turbine $k$ to turbine $l$, $U_{i, k, l}^{m} = 1$.
Otherwise, UAV $i$ needs to fly back to $s$ from the $k$-th turbine and then add another round of inspection starting from turbine $l$.
The detailed procedure is presented in Algorithm \ref{ago:time-cons}.
Quantities $\mathbf{s}_{i, k, l}$, $\theta_{i, k, l}^{s, v}$ and $\mathbf{v}_{i, k, l}$ can be obtained from (\ref{eq:s_veloci_cal}), (\ref{eq:sv_theta}), and (\ref{eq:velocity-relation}), respectively. 
The output of Algorithm \ref{ago:time-cons} is our final routing path.

\begin{algorithm}
\small
\caption{Maximum Flight Time Check}
\label{ago:time-cons}
 \DontPrintSemicolon
\KwIn{$path$, $t_{i}^{max}$}
\KwOut{$U_{i, k, l}^{m}$, $M$}
$m = 1$; $t_{accu} = 0$, $t_{compare} = 0$ \;
\For{ $t = 1$ \KwTo $|{\cal N}_{i}|$ }{
    $k = path[t]$, $l = path[t+1]$ \;
	$t_{accu} = t_{accu} + t_{i, k, l} $\;
    $t_{compare} = t_{accu} + t_{i, l, s} $\;
    \eIf {$t_{accu} < t_{i}^{max}$ \Kwand $t_{compare} < t_{i}^{max}$}{
    $U_{i, k, l}^{m} = 1$ \;
    }
    {
    $U_{i, k, s}^{m} = 1$\;
    $m = m + 1$, $t_{accu} = t_{i, s, l}$, $U_{i, s, l}^{m} = 1$\;  
    }
   }
$M=m$ \;
\end{algorithm}

\subsection{Algorithm Complexity Analysis}

In what follows, the complexity of the proposed algorithms is analyzed.
The computational complexity for the worst-case scenario is provided in our analysis.

First, we analyze the complexity of the algorithms for finding the optimal topology of the UAVs.
Then, a sorting algorithm with the complexity of $n \log (n)$ is utilized.
The complexity of the initialization in Algorithm \ref{ago:deploy_ini} is ${\cal O}(N^{2}) + {\cal O}(NT)$.
Since the number of turbines is the same as the number of UAVs, ${\cal O}(N^{2}) + {\cal O}(NT)$ can be simplified to ${\cal O}(N^{2})$.
On Algorithm \ref{ago:cons_check}, we reduce the number of connections down to $p$.
The complexity of Algorithm \ref{ago:cons_check} is ${\cal O}(N( N  \log(N) + (T-p)))$ considering the worst case is then every UAV is assigned to all turbines.
Line $2$ in Algorithm \ref{ago:UAV_num_mininize} has the complexity of ${\cal O}(N^{2}p + N \log (N))$, where the complexity of the intersection and union is ${\cal O}(p)$.
Each iteration has a complexity of $ {\cal O} (2Np + Np + N^{2}p + N \log(N))$.
Thus, the complexity of Algorithm \ref{ago:UAV_num_mininize} is $ {\cal O} (4N^{2}p + N^{3}p + N^{2} \log(N)+ N \log(N))$.
The proof of the optimality of Algorithms \ref{ago:deploy_ini}-\ref{ago:UAV_num_mininize} is provided in Appendix \ref{subsec:UAV_place_solve_proof}.

Then, we check the complexity of solving the routing problem.
The complexity of executing lines $2$ to $9$ in Algorithm \ref{ago:optimal-route} is ${\cal O}(|{\cal N}_{i}|^{2} \times 2^{|{\cal N}_{i}|})$.
The complexity of lines $11$ to $14$ in Algorithm \ref{ago:optimal-route} is ${\cal O} (|{\cal N}_{i}|)$.
Thus, the complexity of Algorithm \ref{ago:optimal-route} is ${\cal O}(|{\cal N}_{i}|^{2} \times 2^{|{\cal N}_{i}|} + |{\cal N}_{i}|)$.
Algorithm \ref{ago:time-cons} has the complexity ${\cal O} (|{\cal N}_{i}|)$, which depends on the number of turbines assigned to the UAV.
The complexity of applying brute force to find the optimal path is ${\cal O}(|{\cal N}_{i}|\,!)$.
The brute-force method is one of the typical methods used for solving the problem which is NP-hard.
Therefore, the proposed method results in lower complexity compared to brute-force method.
The optimality of Algorithms \ref{ago:optimal-route} and \ref{ago:time-cons} will be proved by comparing with the brute-force method in Section \ref{sec:simulation}. 

\section{Numerical Results}\label{sec:simulation}

\begin{figure}
\begin{center}
\resizebox{3.5in}{!}{%
\includegraphics*{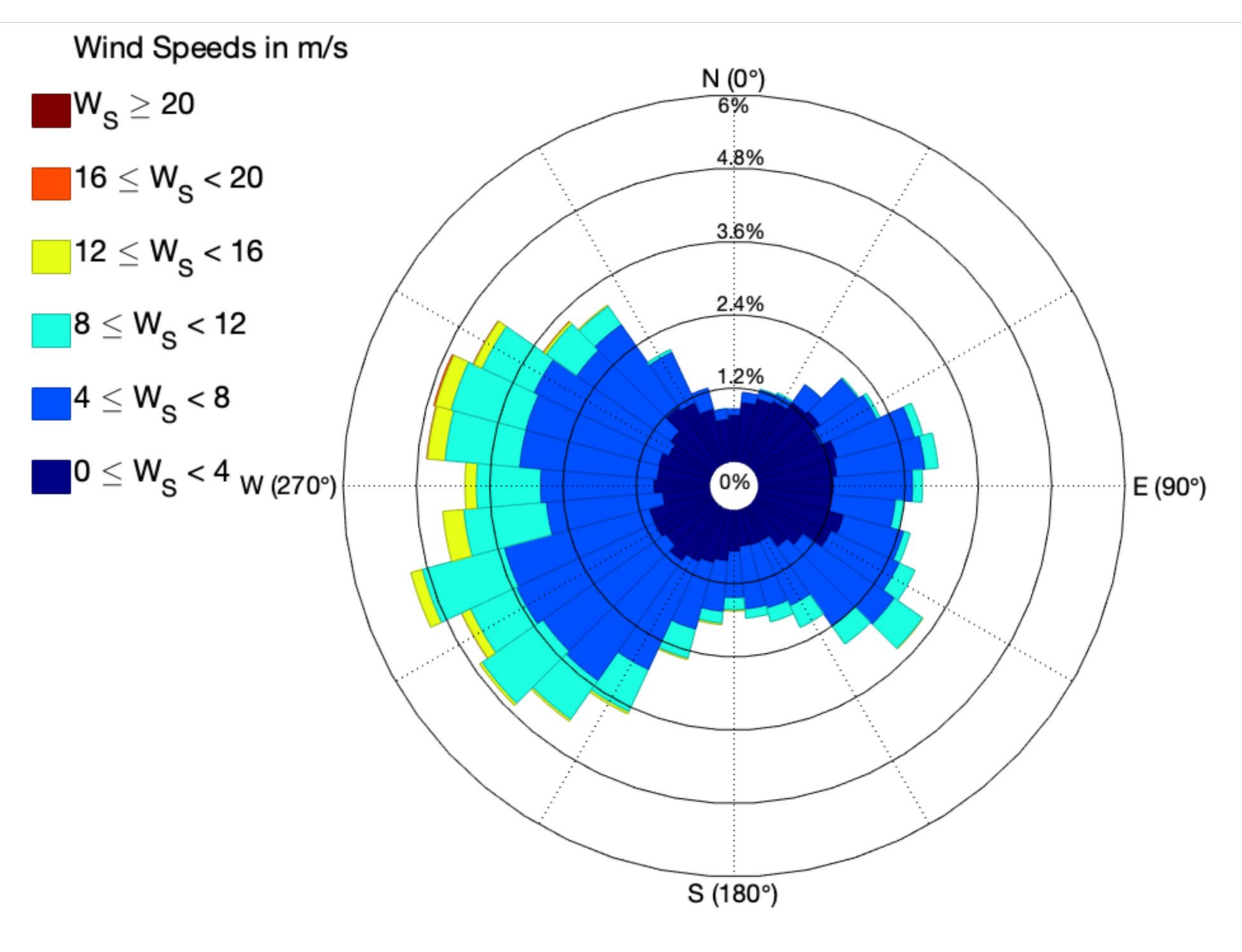} }%
\caption{The wind rose of Walney Island} 
\label{fig:wind_rose}
\end{center}
\end{figure}

In this section, the performance of the proposed method is evaluated based on a real-world dataset.
For the wind farm, we choose the Walney offshore wind farm in the United Kingdom (UK).
The wind farm has an area of $218$ square kilometer (km), a generation capacity of around $1$ Gigawatt (GW), and $189$ turbines.
The data are obtained from Centre for Environmental Data Analysis (CEDA) \cite{Wind_data} and Kingfisher Information Service - Offshore Renewable Cable Awareness (KIS-ORCA) \cite{Turbine_data}.
Specifically, the data from \cite{Wind_data} contain the meteorological measurements in the UK, and then the wind data at Walney Island are utilized in the simulation since wind data inside the wind farm cannot be obtained.
The layout of the Walney wind farm is collected from \cite{Turbine_data}, which contains the longitude and the latitude of each turbine in the wind farm.
The layout of the wind farm obtained from \cite{Turbine_data} is referred to as Walney.
Each turbine is assigned with a code.
The longitude and the latitude of wind turbines are transformed to Cartesian coordinates by using the Mercator projection.
However, a diagram showing the topology of the UAVs with $189$ turbines could lack clarity.
Therefore, we pick $47$ out of $189$ turbines to create another dataset denoted by Walney-1.
The wind data from \cite{Wind_data} can be virtualized with a wind rose as shown in Fig. \ref{fig:wind_rose}.

The UAV used in the simulation is AscTec Falcon 8 \footnote{http://www.asctec.de/en/uav-uas-drones-rpas-roav/asctec-falcon-8/\#pane-0-1}.
This UAV can carry diverse sensors, namely Lidar, ultrasonic sensor, and camera, to inspect turbines.
The maximum speed limit, $u_{i}^{max}$, is set to $16$ m/s.
Then, in the specification of the UAV, it has the maximum flight time between $12$ to $22$ minutes; we set $t_{i}^{max}$ to $20$ minutes.
The model introduced in Appendix \ref{subsec:UAV_energy_model} can be used to verify these settings.
The maximum resistance to the wind speed of the UAV is $u_{i}^{wind} = 15$ m/s.
The UAVs need to communicate in order to prevent collision.
Therefore, the maximum distance of the communication between UAVs is set to $5$ km.
Every UAV can be assigned to inspect up to $5$ turbines, i.e., $p = 5$.
The range between $0$ and $2\pi$ is discretized into $\mu = 36$ segments.

\subsection{The Topology of the UAVs in the Wind Farm} \label{subsec:result_topology}

\begin{figure}
\begin{center}
\resizebox{3in}{!}{%
\includegraphics*{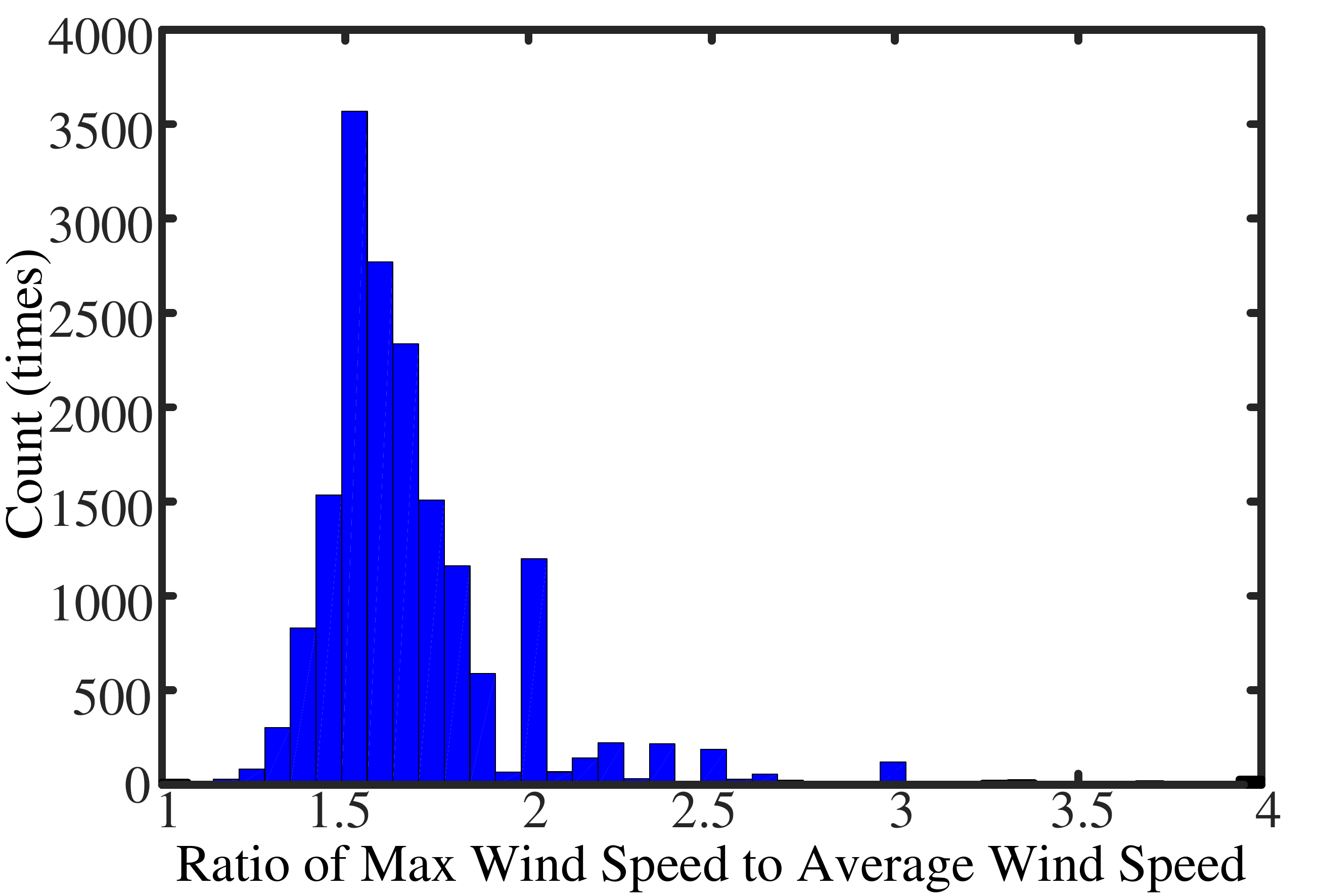} }%
\caption{The ratio of hourly wind gust to hourly average wind speed }
\label{fig:wind_max_ave}
\end{center}
\end{figure}

To solve the placement problem, we need to determine the flying range and then use it as input.
The histogram in Fig. \ref{fig:wind_max_ave} counts the occurrences of the quotient of hourly wind gust and hourly average wind speed.
Values below $2$ account for $93.14 \%$ of the data. 
Thus, there is a higher than $90 \%$ chance for the peak value not to exceed twice the average value.
Given a maximum wind speed resistance of up to $15$ m/s for this UAV, it is plausible to set the hourly average wind speed up to which the UAV is allowed to perform inspection to $\epsilon_{v} = 8$ m/s.
Wind speeds above $15$ m/s are rarely encountered.
With $\epsilon_{v}$, Algorithm \ref{ago:range_deter} is applied to determine the flying range.

\begin{figure*}[t!] 
\subfloat[$step = 0$\label{fig:step0}]
{\scalebox{0.19}{\includegraphics{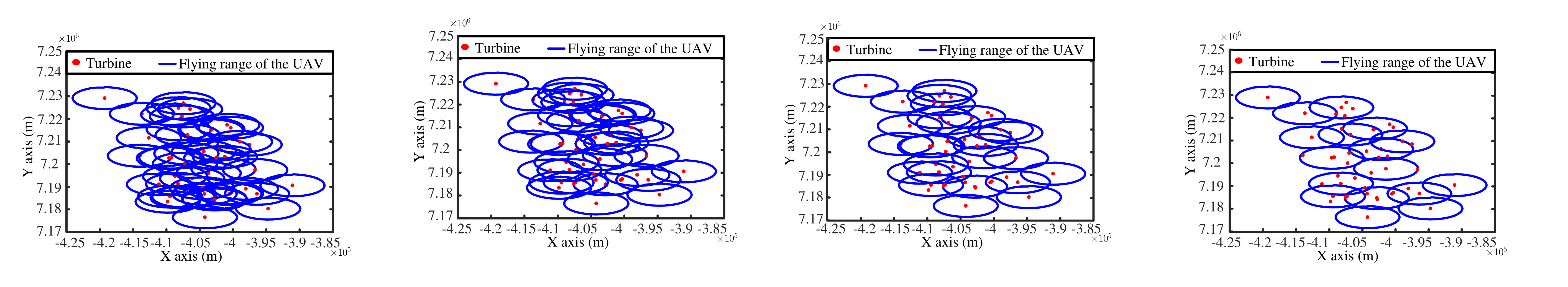}}}
\hspace{-0.39cm}
\subfloat[$step = 10$\label{fig:step10}]
{\scalebox{0.19}{\includegraphics{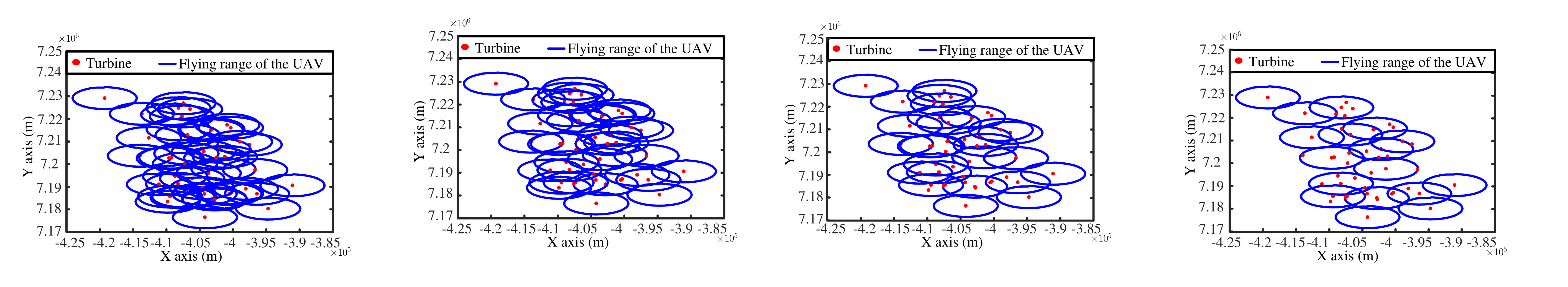}}}
\hspace{-0.39cm}
\subfloat[$step = 20$\label{fig:step20}]
{\scalebox{0.19}{\includegraphics{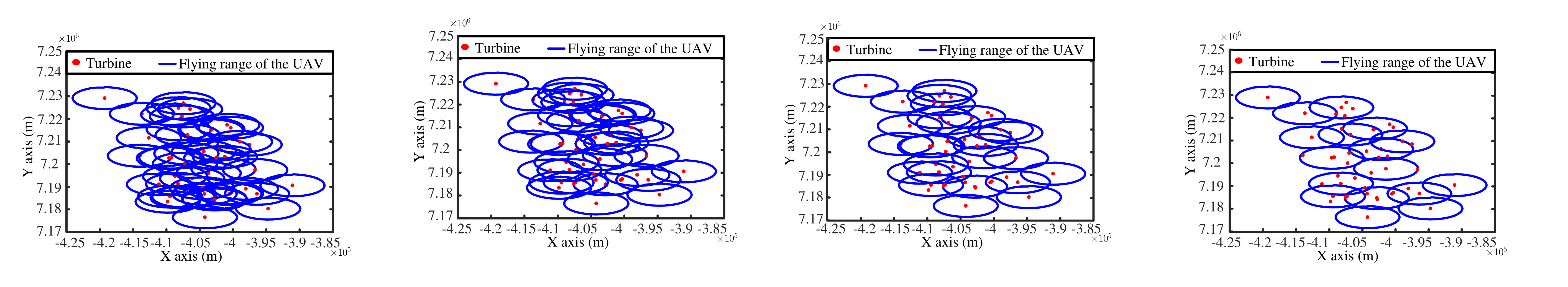}}}
\hspace{-0.43cm}
\subfloat[$step = 30$\label{fig:step30}]
{\scalebox{0.19}{\includegraphics{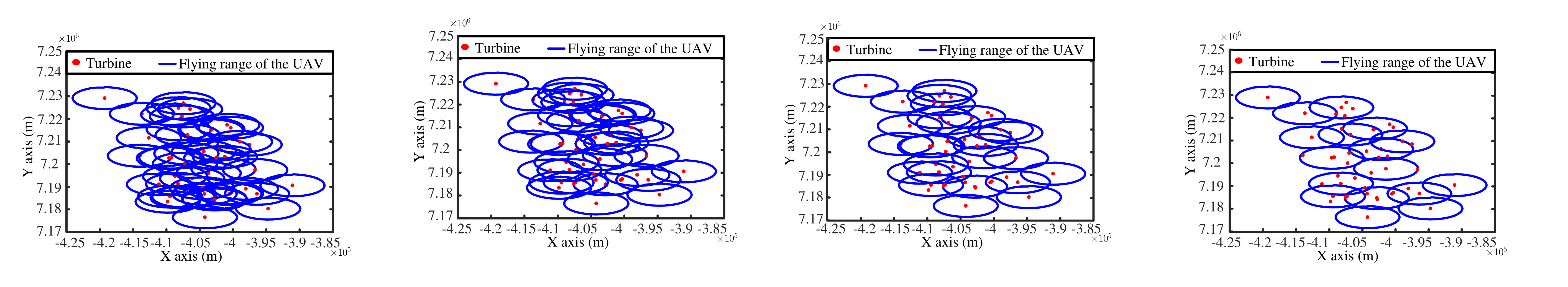}}}
\caption{UAV placement process.}  \label{fig:deploy_result}
\end{figure*}

The flying range is then applied to determine the topology.
Fig. \ref{fig:deploy_result} shows the iteration process and the final placement results by applying the proposed algorithms to the Walney-1 dataset.
At the beginning of the iteration, $47$ UAVs are placed in the wind farm as shown in Fig. \ref{fig:step0}.
In Fig. \ref{fig:step10} and \ref{fig:step20}, redundant UAVs are deleted based on the proposed algorithms.
After $30$ iterations, the algorithm stops and outputs the final placement results as shown in Fig. \ref{fig:step30}.
According to our results, only $17$ UAVs are required to cover all wind turbines in the offshore wind farm.
We also observe that all the turbines are assigned to a UAV, and every UAV serves no more than $5$ turbines.
The same setting are applied to solve the placement problem with the Walney dataset.
The results reveal that we need to place $63$ UAVs to cover all turbines in the wind farm.

\begin{table} \small
\begin{center}
\caption{The placement results of changing the number a UAV can serve}
\label{tb:change_num_turbine}
\begin{tabular}{|c|c|c|c|c|c|c|c|}
\hline
Dataset   &  $T$ & $t_{i}^{max}$  & $\epsilon_{v}$ (m/s) & $p$ &  $\sum_{i}^{N} a_{i}$   \\
\hline
\multirow{6}{*}{Walney-1}   & \multirow{3}{*}{$47$}  & \multirow{3}{*}{$20$} &  \multirow{3}{*}{$6$}    & $5$   & $16$           \\
\cline{5-6}      &  & & & $8$ & $15$ \\
\cline{5-6}      &  & & & $10$ & $15$ \\
\cline{2-6}
& \multirow{3}{*}{$47$}  & \multirow{3}{*}{$20$} &  \multirow{3}{*}{$8$}    & $5$   & $17$           \\
\cline{5-6}      &  & & & $8$ & $16$ \\
\cline{5-6}      &  & & & $10$ & $16$ \\
\hline \hline
\multirow{6}{*}{Walney}  & \multirow{3}{*}{189}  & \multirow{3}{*}{$20$} &\multirow{3}{*}{$6$}  & $5$      & $61$    \\
\cline{5-6}    & &  &  & $8$ & $45$ \\
\cline{5-6}    & &  &  & $10$ & $42$ \\
\cline{2-6}
 & \multirow{3}{*}{189}  & \multirow{3}{*}{$20$} &\multirow{3}{*}{$8$}  & $5$      & $63$    \\
\cline{5-6}    & &  &  & $8$ & $56$ \\
\cline{5-6}    & &  &  & $10$ & $54$ \\
\hline

\end{tabular}
\end{center}
\end{table}

In the proposed algorithms, the number of the placed UAVs can be further reduced if one UAV can serve more turbines by increasing the value of $p$.
On the other hand, increasing the flying range by decreasing the value of $\epsilon_{v}$ can also reduce the number of the deployed UAVs.
Therefore, we compare the influence of $p$ and $\epsilon_{v}$ to the placement results in Table \ref{tb:change_num_turbine}.
For Walney-1, there is no significant difference because this dataset includes only $47$ turbines.
However, for Walney itself, $63$ UAVs are required according to the previous setting, $p=5$ and $\epsilon_{v} = 8$.
The number of placed UAVs can be cut down to $42$ if we change $p$ and $\epsilon_{v}$ to $10$ and $6$ m/s, respectively. 
However, with more turbines to be served per UAV, it will be more challenging to find an optimal routing solution.
This issue will be discussed in the next section.

\subsection{Routing Result}

\begin{figure}
\begin{center}
\resizebox{3in}{!}{%
\includegraphics*{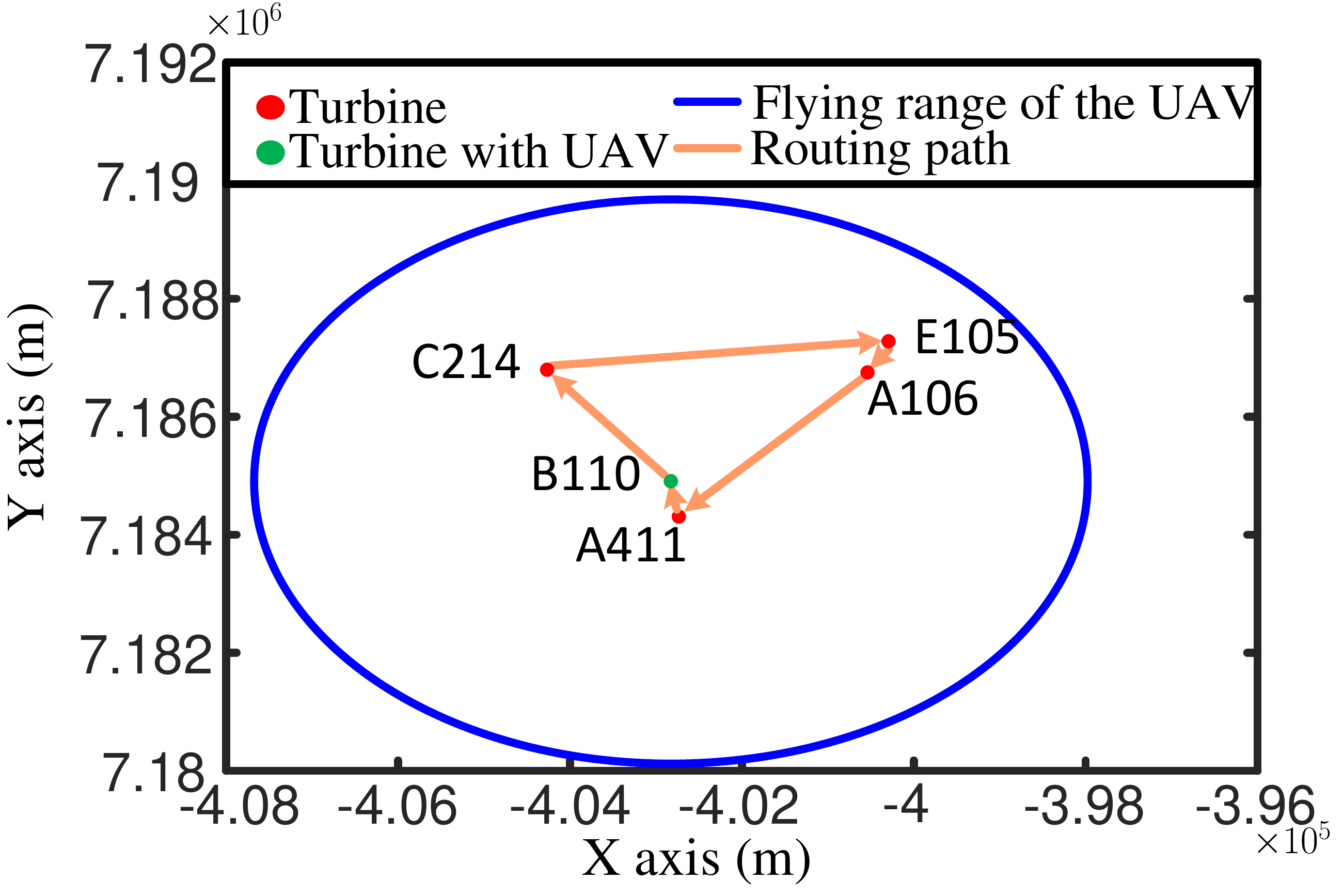} }%
\caption{The routing path of UAV $15$} 
\label{fig:Routing_UAV_2}
\end{center}
\end{figure}

The results of the placement is applied to show how to route the UAVs to inspect the wind turbines.
In this case, $w_{s}$ is set to $8$ m/s and $\theta_{w}^{met}$ is set to $\frac{\pi}{2}$ (east wind).
We take UAV $15$ as an example.
The UAV is placed at the turbine whose code is B110.
The result is shown in Fig. \ref{fig:Routing_UAV_2}.
To minimize inspection time, the UAV should avoid facing the head wind. 
Therefore, the UAV goes to C214 first and then chooses E105 afterwards.
After E105, the UAV uses the tail wind to go to A106 and A411.
The total flight time for the inspection is $15.30$ minutes.
The proposed method is compared with the brute-force method and the branch-and-bound method, which are common algorithms to solve MILP problem, in Table \ref{tb:compare_diff_sol_rout}.
According to the results, Algorithms \ref{ago:optimal-route} and \ref{ago:time-cons} obtain the same result as the brute-force method.
However, this is not true for the branch-and-bound method.
In this case, the optimality of Algorithms \ref{ago:optimal-route} and \ref{ago:time-cons} can be proved.

\begin{table} \small
\begin{center}
\caption{Total flight time under different solutions}
\label{tb:compare_diff_sol_rout}
\renewcommand{\arraystretch}{1.1}
\begin{tabular}{|c|c|c|c|c|c|c|c|}
\hline
\multirow{2}{*}{Method}    & \multirow{2}{*}{\shortstack{Total flight\\ time (mins)}}   &  \multirow{2}{*}{$path$ }  \\
& & \\
\hline
\multirow{2}{*}{Branch-and-bound}               & \multirow{2}{*}{$15.7343$}   & \multirow{2}{*}{ \shortstack[l]{B110$>$A106$>$E105$>$\\C214$>$A411$>$B110}}   \\
& & \\
\hline 
\multirow{2}{*}{Brute-force}                  & \multirow{2}{*}{$15.3049$}   & \multirow{2}{*}{ \shortstack[l]{B110$>$C214$>$E105$>$\\A106$>$A411$>$B110}}   \\
& & \\
\hline 
\multirow{2}{*}{Algorithm \ref{ago:optimal-route} and \ref{ago:time-cons}}              & \multirow{2}{*}{$15.3049$}   & \multirow{2}{*}{ \shortstack[l]{B110$>$C214$>$E105$>$\\A106$>$A411$>$B110}}    \\
& & \\
\hline 
\end{tabular}
\end{center}
\end{table}

Fig. \ref{fig:Routing_UAV_2} and Table \ref{tb:change_max_time} show that the UAV can finish the inspection in one round.
This is because $t_{i}^{max}$ is set to $20$ minutes.
Of interest is also the performance of Algorithm \ref{ago:time-cons} when reducing $t_{i}^{max}$.
For that, we decrease the maximum flight time to $15$ and $12$ minutes.
The results are recorded in Table \ref{tb:change_max_time}.
The UAV needs two rounds if $t_{i}^{max}$ is set to $15$ minutes.
For the first round, it takes $14.44$ minutes; it is $1.27$ minutes for the second round.
The UAV still needs two rounds to finish the inspection if $t_{i}^{max}$ is further reduced to $12$ minutes.
The first and the second rounds take $6.54$ and $10.73$ minutes, respectively.
However, the turbines assigned to each route vary in each round compared to $t_{i}^{max} = 15$ minutes as shown in Table \ref{tb:change_max_time}.
By applying the proposed algorithm, we can ensure that the UAV does not spend more than $t_{i}^{max}$ on inspection.

\begin{table} \small
\begin{center}
\caption{The routing results of changing the maximum flight time}
\label{tb:change_max_time}
\begin{tabular}{|c|c|c|c|c|c|c|c|}
\hline
$t_{i}^{max}$ (mins)    & $m$   &  $path$   \\
\hline
$20$                    & $1$   & B110$>$C214$>$E105$>$A106$>$A411$>$B110           \\
\hline \hline
\multirow{2}{*}{$15$}   & $1$ & B110$>$C214$>$E105$>$A106$>$B110     \\
\cline{2-3}
        & $2$   & B110$>$A411$>$B110 \\
\hline
\hline
\multirow{2}{*}{$12$}   & $1$ & B110$>$C214$>$B110     \\
\cline{2-3}
        & $2$   & B110$>$E105$>$A106$>$A411$>$B110 \\
\hline
\end{tabular}
\end{center}
\end{table}

As mentioned in Section \ref{subsec:result_topology}, increasing $p$ may cause an issue for routing.
This issue is discussed here, and $p$ is set to $10$.
Two UAVs, with $i=15$ and $i=5$, are taken out to compare the required number of routes in Table \ref{tb:change_num_turbine_time}.
UAV $15$ serves $5$ and UAV $5$ serves $9$ turbines.
We consider $t_{i}^{max}$ to be $15$ and $20$ minutes.
Then, $w_{s}$ is set to $8$ m/s and $\theta_{w}^{met}$ to $0$ and $\pi$.
For UAV $15$, it only needs one round to finish the inspection under four different settings.
However, UAV $5$ needs $2$ rounds if $t_{i}^{max}$ is $20$ minutes.
The number of the required routes doubles to $4$ if $t_{i}^{max}$ is reduced to $15$ minutes.
Therefore, the UAV needs to spend more time and number of routes on routing if one UAV is to inspect more turbines.
Also, in the wind data, $83.02\%$ of the hourly average wind speed is lower than $8$ m/s ($\epsilon_{v} = 8$). 
However, the value is $63.06 \%$ for the hourly average wind speed lower than $6$ m/s ($\epsilon_{v} = 6$).
Therefore, if we reduce $\epsilon_{v}$, it means the UAVs have less chance to start the inspection.
From these results, it is clear that there exists a fundamental tradeoff between maximum flight time, number of turbines one UAV can serve, and the flying range.
In the case of Walney wind farm, $p=5$ and $\epsilon_{v} = 8$ m/s is better than $p=10$ and $\epsilon_{v} = 6$ m/s.

\begin{table} \small
\begin{center}
\caption{The routing results of changing $t_{i}^{max}$ and $p$}
\label{tb:change_num_turbine_time}
\begin{tabular}{|c|c|c|c|c|c|c|c|}
\hline
$i$   &  $|{\cal N}_{i}|$  & $p$ & $\epsilon_{v}$ & $t_{i}^{max} (min)$  & $w_{s} (m/s)$ & $\theta_{w}^{met}$ &$M$   \\
\hline
\multirow{4}{*}{$15$} & \multirow{4}{*}{$5$}  & \multirow{4}{*}{$5$} & \multirow{4}{*}{$8$} &  \multirow{2}{*}{$15$}    & $8$   & $0$       &  $1$  \\
\cline{6-8}     & &  & & & $8$ & $\pi$ & $1$ \\
\cline{5-8} 
& & &  & \multirow{2}{*}{$20$}    & $8$   & $0$     & $1$ \\
\cline{6-8}
& & & & & $8$ & $\pi$ & $1$\\
\hline \hline
\multirow{4}{*}{$5$} & \multirow{4}{*}{$9$}  & \multirow{4}{*}{$10$} & \multirow{4}{*}{$8$} & \multirow{2}{*}{$15$}    & $8$   & $0$       &   $4$  \\
\cline{6-8}    &  &  & & & $8$ & $\pi$ &  $4$ \\
\cline{5-8} 
& & & & \multirow{2}{*}{$20$}    & $8$   & $0$       &  $2$ \\
\cline{6-8}
& & & & & $8$ & $\pi$ & $2$ \\
\hline 
\end{tabular}
\end{center}
\end{table}

\section{Conclusion}\label{sec:conclusion}

In this paper, we presented a framework for utilizing UAVs to inspect the wind turbines in an offshore wind farm.
An optimization problem was formulated to minimize the number of UAVs to be placed in the offshore wind farm with consideration of the challenging offshore wind condition faced by the UAVs.
Another formulated optimization problem was to find an optimal route for the UAVs to inspect the wind turbines.
We designed heuristic algorithms to solve both problems and analyzed the complexity of the proposed algorithms.
For the purpose of validation, real-world data were utilized (meteorological measurements recorded by Centre for Environmental Data Analysis (CEDA), and positions of the turbines through Kingfisher Information Service - Offshore Renewable Cable Awareness (KIS-ORCA)).
With the proposed methods, we can discover how many UAVs are needed to automatically inspect the turbines in an offshore wind farm.
The optimal routing path can also be obtained for the inspection under different wind conditions.
With the proposed framework, more efficient and more frequent inspection of wind turbines can be achieved for the wind farm operators.
By doing so, the loss due to failures of the wind turbines can be reduced.
In our future work, we will study how UAVs can bring more benefits to the operation of offshore wind farms.


\appendix

\subsection{Flying Distance Determination} \label{subsec:UAV_energy_model}

In (\ref{eq:flying_range_max}), the flying distance is mainly determined by the product of $u_{i}^{max}$ and $t_{i}^{max}$.
However, $u_{i}^{max}$ and $t_{i}^{max}$ may influence each other.
Therefore, for $u_{i}^{max}$, it is reasonable to be assigned with the maximum flying speed listed in the specification of the UAV.
Quantity $t_{i}^{max}$ should be calculated based on $u_{i}^{max}$, and it can be obtained from calculating the energy consumption of the UAV.
Before introducing the model, $V$ is used to represent the airspeed, $||\mathbf{v}_{i, k, l}||_{2}$, for the sake of notational simplicity.
The power consumption of a UAV flying with airspeed $V$ can then be modeled \cite{2001-UAV-energy-model,2006-UAV-energy-model} as
\begin{equation}\label{eq:UAV_power_model}
\begin{array}{ll}
P(V) & = \underbrace{P_{o} \left( 1 + \frac{ 3 V^{2} }{ U_{tip}^{2} }  \right)}_{ \mbox{blade profile} } + \underbrace{P_{i} \left( \sqrt{ 1 + \frac{V^{4}}{4v_{o}^{4}}}  - \frac{V^{2}}{2v_{o}^{2}} \right)}_{ \mbox{induced power} }  \\
& + \underbrace{\frac{1}{2} d_{0} sol \rho A_{disc} V^{3}}_{ \mbox{parasite} },
\end{array}
\end{equation}
where $P_{0}$ and $P_{i}$ are two constants defined in (\ref{eq:p0_pi_equation}) representing the blade profile power and induced power in hovering status, respectively.
\begin{equation} \label{eq:p0_pi_equation}
\left\{
\begin{array}{ll}
                P_{0} = \frac{\delta}{8} sol \rho A_{disc} \Omega^{3} R^{3},    \\
                P_{i} = (1+k_{cor})\frac{W^{3/2}}{\sqrt{2\rho A_{disc}}} .
\end{array}
        \right.
\end{equation}
Quantity $U_{tip}$ denotes the tip speed of the rotor blade, and $v_{0}$ is known as the mean rotor induced velocity in hover.
The fuselage drag ratio and rotor solidity are denoted by $d_{0}$ and $sol$, respectively.
The air density is denoted by $\rho$, and $A_{disc}$ is the rotor disc area.
The maximum flight time can then be obtained from 
\begin{equation} \label{eq:ti_max_cal}
t_{i}^{max} =  \frac{P_{bat}}{P(V)}, 
\end{equation}
where $P_{bat}$ is the battery capacity.
The parameters and the meanings used in the calculation are provided in Table \ref{tb:max_fly_time_cal}.
Moreover, the values of these parameters are based on the values on the specification of the UAV that is used in the simulation.
Based on (\ref{eq:UAV_power_model}), (\ref{eq:p0_pi_equation}), and Table \ref{tb:max_fly_time_cal}, the energy consumption of the UAV is $212.82$ Watt.
Then, according to (\ref{eq:ti_max_cal}), $t_{i}^{max}$ is $20.02$ minutes.
Therefore, $t_{i}^{max}$ is set to $20$ minutes in Section \ref{sec:simulation}.

\begin{table} \small
\renewcommand{\arraystretch}{1.2}
\begin{center}
\caption{Parameters and their meanings for calculating maximum flight time}
\label{tb:max_fly_time_cal}
\begin{tabular}{|c|l|c|c|c|c|c|c|}
\hline
Parameter   & \multicolumn{1}{c|}{Physical Meaning}   &  Value  \\
\hline
$V$    & Air speed in in meter per second (m/s)  &  $16$   \\
\hline
$\omega$    & UAV weight in Newton   &  $16$   \\
\hline
$R$    & Rotor radius in m    & $0.1016$ \\
\hline
$\rho$    & Air density in $\mbox{kg}/\mbox{m}^{3}$  & $1.2250$ \\
\hline
$A_{disc}$    & Rotor disc area in $\mbox{m}^{2}$, $A_{disc} \triangleq \pi R^{2}$  & $0.0314$ \\
\hline
\multirow{2}{*}{$\Omega$}    &  \multirow{2}{*}{\shortstack[l]{Angular velocity of UAV blade in radian \\ per second (rad/s) } } & \multirow{2}{*}{$300$}\\
&  & \\
\hline
$U_{tip}$    & Tip speed of the rotor blade, $U_{tip} \triangleq \Omega R$   & $30$ \\
\hline
$b_{num}$    & Number of blade   & $8$ \\
\hline
$cord$    & Chord length of UAV blade in m    & $0.09$ \\
\hline
$sol$    & Rotor solidity, $sol \triangleq  \frac{b_{num}cord}{\pi R}$   & $2.5464$ \\
\hline
\multirow{2}{*}{$k_{cor}$ }   & \multirow{2}{*}{ \shortstack[l]{ Incremental correction factor to\\induced power} }   & \multirow{2}{*}{$0.1$} \\
&  & \\
\hline
\multirow{2}{*}{$v_{0}$ }   & \multirow{2}{*}{ \shortstack[l]{Mean rotor induced velocity in hover,\\ $v_{0} \triangleq \sqrt{ \omega/(2 \rho A_{disc}) } $} }   & \multirow{2}{*}{$14.4179$} \\
&  & \\
\hline
$\delta$    & Profile drag coefficient   & $0.0120$ \\
\hline
$S_{fp}$    & Fuselage equivalent flat plate area in $\mbox{m}^{2}$   & $0.0063$ \\
\hline
$d_{0}$    & Fuselage drag ratio, $d_{0} \triangleq \frac{S_{fp}}{A sol} $   & $0.0787$ \\
\hline
\multirow{2}{*}{$P_{bat}$ }   & \multirow{2}{*}{ \shortstack[l]{ Battery capacity of the UAV in Ampere\\ hour (Ah)} }   & \multirow{2}{*}{$6.25$} \\
&  & \\
\hline

\end{tabular}
\end{center}
\end{table}

\subsection{Proof of Optimality of Algorithm 2-4 for Solving Placement Problem} \label{subsec:UAV_place_solve_proof}

Suppose the optimal solution of the placement problem is $\mathbf{A}^{*}=[a_{i}^{*}]_{1\times N}$, $\mathbf{B}^{*} = [b_{i, k}^{*}]_{N \times T}$, $\mathbf{C}^{*} = [c_{i, j}^{*}]_{N \times N}$, and $[x_{i}^{*}, y_{i}^{*}]$ for $i \in \{i| a_{i}^{*}=1\}$.
Also, $n^{*} = \sum_{i=1}^{N} a_{i}^{*}$ is the optimal number of UAVs that should be deployed in the wind farm.
On the other hand, the solution obtained from the proposed algorithms is denoted by $\mathbf{A}$, $\mathbf{B}$, $\mathbf{C}$, and $[x_{i}, y_{i}]$ for $i \in \{i| a_{i}=1\}$.
Then, $n$ is defined by $ \sum_{i=1}^{N} a_{i}$.
With these notations, this implies that $n > n^{*}$ and at least one UAV can be removed from the solution obtained from the proposed algorithms.
That is, UAV $j$ can be removed if constraints in (\ref{eq:UAV_deploy_problem}) are still satisfied without UAV $j$.
Thus, it is necessary to recheck if there is still an overlap between ${\cal N}_{i}$ and ${\cal N}_{j}$.
However, the end condition of Algorithm $4$ is that there is no overlap between ${\cal N}_{i}$ and ${\cal N}_{j}$.
This contradiction implies that $n$ is the same as $n^{*}$.
Also, $\mathbf{A} = \mathbf{A}^{*}$, $\mathbf{B} = \mathbf{B}^{*}$, $\mathbf{C} = \mathbf{C}^{*}$, and $[x_{i}, y_{i}] = [x_{i}^{*}, y_{i}^{*}]$.
Finally, we can conclude that the algorithm would converge to the optimal solution.

{\renewcommand{\baselinestretch}{1}
\begin{footnotesize}
\bibliographystyle{IEEEtran}
\bibliography{References_TII_UAV}

\begin{thebibliography}{10}
\providecommand{\url}[1]{#1}
\csname url@samestyle\endcsname
\providecommand{\newblock}{\relax}
\providecommand{\bibinfo}[2]{#2}
\providecommand{\BIBentrySTDinterwordspacing}{\spaceskip=0pt\relax}
\providecommand{\BIBentryALTinterwordstretchfactor}{4}
\providecommand{\BIBentryALTinterwordspacing}{\spaceskip=\fontdimen2\font plus
\BIBentryALTinterwordstretchfactor\fontdimen3\font minus
  \fontdimen4\font\relax}
\providecommand{\BIBforeignlanguage}[2]{{%
\expandafter\ifx\csname l@#1\endcsname\relax
\typeout{** WARNING: IEEEtran.bst: No hyphenation pattern has been}%
\typeout{** loaded for the language `#1'. Using the pattern for}%
\typeout{** the default language instead.}%
\else
\language=\csname l@#1\endcsname
\fi
#2}}
\providecommand{\BIBdecl}{\relax}
\BIBdecl

\bibitem{wind-power-forecast}
\BIBentryALTinterwordspacing
{Wood Mackenzie}, ``{Global wind power market outlook update: Q2 2019},'' Wood
  Mackenzie Power \& Renewables, Tech. Rep., Jun. 2019. [Online]. Available:
  \url{https://www.woodmac.com/reports/power-markets-global-wind-power-market-outlook-update-q2-2019-318297/}
\BIBentrySTDinterwordspacing

\bibitem{wind-power-europe}
\BIBentryALTinterwordspacing
{Wind Europe}, ``{Europe Installs 4.9 GW of New Wind Energy Capacity in First
  Half of 2019},'' Tech. Rep., Jul. 2019. [Online]. Available:
  \url{https://windeurope.org/newsroom/press-releases/europe-installs-4-9-gw-of-new-wind-energy-capacity-in-first-half-of-2019/}
\BIBentrySTDinterwordspacing

\bibitem{2007-turbine-failure}
{J. Ribrant and L. M. Bertling}, ``{Survey of Failures in Wind Power Systems
  With Focus on Swedish Wind Power Plants During 1997–2005},'' \emph{IEEE
  Trans. Energy Convers.}, vol.~22, no.~1, pp. 167--173, Mar. 2007.

\bibitem{2014-DTU-report}
{K. Branner and A. Ghadirian}, ``{Database about Blade Faults},'' DTU Wind
  Energy E-0067, Tech. Rep. 978-87-93278-09-7, Dec. 2014.

\bibitem{2015-fibre-lee}
{J.-K. Lee, J.-Y. Park, K.-Y. Oh, S.-H. Ju, and J.-S. Lee}, ``{Transformation
  Algorithm of Wind Turbine Blade Moment Signals for Blade Condition
  Monitoring},'' \emph{Renew. Energy}, vol.~79, pp. 209--218, Jul. 2015.

\bibitem{2016-lidar-inspect-UAV}
{B. E. Sch{\"a}fer, D. Picchi, T. Engelhardt, and D. Abel}, ``{Multicopter
  Unmanned Aerial Vehicle for Automated Inspection of Wind Turbin},'' in
  \emph{Proc. Medit. Conf. Control Autom. (MED)}, Athens, Greece, Jun. 2016,
  pp. 244--249.

\bibitem{2018-blade-thermal-yang}
{R. Yang, Y. He, A. Mandelis, N. Wang, X. Wu, and S. Huang}, ``{Induction
  Infrared Thermography and Thermal-Wave-Radar Analysis for Imaging Inspection
  and Diagnosis of Blade Composites},'' \emph{IEEE Trans. Ind. Informat.},
  vol.~14, no.~12, pp. 5637--5647, Dec. 2018.

\bibitem{2017-blade-mmwave-wang}
{J. R. Gallion and R. Zoughi}, ``{Millimeter-Wave Imaging of Surface-Breaking
  Cracks in Steel With Severe Surface Corrosion},'' \emph{IEEE Trans. Instrum.
  Meas.}, vol.~66, no.~10, pp. 2789--2791, Oct. 2017.

\bibitem{2017-gearbox-failure-DNN}
{L. Wang, Z. Zhang, H. Long, J. Xu, and R. Liu}, ``{Wind Turbine Gearbox
  Failure Identification With Deep Neural Networks},'' \emph{IEEE Trans. Ind.
  Informat.}, vol.~13, no.~3, pp. 1360--1368, Jun. 2017.

\bibitem{2019-blade-ice}
\BIBentryALTinterwordspacing
{B. Yuan, C. Wang, F. Jiang, M. Long, P. S. Yu, and Y. Liu}, ``{WaveletFCNN: A
  Deep Time Series Classification Model for Wind Turbine Blade Icing
  Detection},'' \emph{arXiv}, Feb. 2019. [Online]. Available:
  \url{https://arxiv.org/abs/1902.05625}
\BIBentrySTDinterwordspacing

\bibitem{2015-UAV-PV-panel}
{M. Aghaei, F. Grimaccia, C. A. Gonano, and S. Leva}, ``{Innovative Automated
  Control System for PV Fields Inspection and Remote Control},'' \emph{IEEE
  Trans. Ind. Electron.}, vol.~62, no.~11, pp. 7287--7296, Nov. 2015.

\bibitem{2017-UAV-blade-wang}
{L. Wang and Z. Zhang}, ``{Automatic Detection of Wind Turbine Blade Surface
  Cracks Based on UAV-Taken Images},'' \emph{IEEE Trans. Ind. Electron.},
  vol.~64, no.~9, pp. 7293--7303, Sep. 2017.

\bibitem{2017-power-reading-UAV}
{J. R. T. Neto, A. Boukerche, R. S. Yokoyama, D. L. Guidoni, R. I. Meneguette,
  J. Ueyama, and L. A. Villas}, ``{Performance Evaluation of Unmanned Aerial
  Vehicles in Automatic Power Meter Readings},'' in \emph{Proc. Ad Hoc New.},
  vol.~60, no.~C, Amsterdam, The Netherlands, May 2017, pp. 11--25.

\bibitem{2019-UAV-blade-wang}
{L. Wang, Z. Zhang, and X. Luo}, ``{A Two-Stage Data-Driven Approach for
  Image-Based Wind Turbine Blade Crack Inspections},'' \emph{IEEE/ASME Trans.
  Mechatronics}, vol.~24, no.~3, pp. 1271--1281, Jun. 2019.

\bibitem{2018-UAV-power-inspect}
{G. J. Lim, S. Kim, J. Cho, Y. Gong, and A. Khodaei}, ``{Multi-UAV
  Pre-Positioning and Routing for Power Network Damage Assessment},''
  \emph{IEEE Trans. Smart Grid}, vol.~9, no.~4, pp. 3643--3651, Jul. 2018.

\bibitem{2020-UAV-insulator-inspect}
{X. Tao, D. Zhang, Z. Wang, X. Liu, H. Zhang, and D. Xu}, ``{Detection of Power
  Line Insulator Defects Using Aerial Images Analyzed With Convolutional Neural
  Networks},'' \emph{IEEE Trans. Syst., Man, Cybern. Syst.}, vol.~50, no.~4,
  pp. 1486--1498, Apr. 2020.

\bibitem{2019-UAV-power-inspect}
{S. Y. Derakhshandeh, Z. Mobini, M. Mohammadi, and M. Nikbakht},
  ``{UAV-Assisted Fault Location in Power Distribution Systems: An Optimization
  Approach},'' \emph{IEEE Trans. Smart Grid}, vol.~10, no.~4, pp. 4628--4636,
  Jul. 2019.

\bibitem{2018-UAV-rui-tsg}
{J. Zhang, Y. Zeng, and R. Zhang}, ``{UAV-Enabled Radio Access Network:
  Multi-Mode Communication and Trajectory Design},'' \emph{IEEE Trans. Sig.
  Proc.}, vol.~66, no.~20, pp. 5269--5284, Oct. 2018.

\bibitem{2019-UAV-rui-twc-secure}
{Y. Zeng, J. Xu, and R. Zhang}, ``{Securing UAV Communications via Joint
  Trajectory and Power Control},'' \emph{IEEE Trans. Wireless Commun.},
  vol.~18, no.~2, pp. 1376--1389, Feb. 2019.

\bibitem{2019-UAV-rui-twc-energy}
------, ``{Energy Minimization for Wireless Communication With Rotary-Wing
  UAV},'' \emph{IEEE Trans. Wireless Commun.}, vol.~18, no.~4, pp. 2329--2345,
  Apr. 2019.

\bibitem{2019-UAV-deploy-zhang}
{X. Zhang and L. Duan}, ``{Fast Deployment of UAV Networks for Optimal Wireless
  Coverage},'' \emph{IEEE Trans. Mobile Comput.}, vol.~18, no.~3, pp. 588--601,
  Mar. 2019.

\bibitem{2019-UAV-deploy-sun}
{J. Sun and C. Masouros}, ``{Deployment Strategies of Multiple Aerial BSs for
  User Coverage and Power Efficiency Maximization},'' \emph{IEEE Trans.
  Commun.}, vol.~67, no.~4, pp. 2981--2994, Apr. 2019.

\bibitem{2018-deplpoy-UAV-predictive}
\BIBentryALTinterwordspacing
{Q. Zhang, W. Saad, M. Bennis, X. Lu, M. Debbah, and W. Zuo}, ``{Predictive
  Deployment of UAV Base Stations in Wireless Networks: Machine Learning Meets
  Contract Theory},'' \emph{arXiv}, Nov. 2018. [Online]. Available:
  \url{https://arxiv.org/abs/1811.01149}
\BIBentrySTDinterwordspacing

\bibitem{2016-UAV-dockstation}
\BIBentryALTinterwordspacing
L.~Mathews, ``{Airobotics Drone Recharging Station will Function Forever
  without Humans},'' Jun. 2016. [Online]. Available:
  \url{https://www.geek.com/tech/airobotics-drone-recharging-station-will-function-forever-without-humans-1658919/}
\BIBentrySTDinterwordspacing

\bibitem{2007-UAV-Adhoc-collision}
{X. Wang, V. Yadav, and S. N. Balakrishnan}, ``{Emission-Aware and
  Cost-Effective Distributed Demand Response System for Extensively Electrified
  Large Ports},'' \emph{IEEE Trans. Control Syst. Technol.}, vol.~15, no.~4,
  pp. 672--679, Jul. 2007.

\bibitem{2018-UAV-Adhoc-collision}
{F. Fabra, C. T. Calafate, J. C. Cano, and P. Manzoni}, ``{A Collision
  Avoidance Solution for UAVs Fllowing Planned Missions},'' in \emph{Proc. IEEE
  Wireless Commun. Netw. Conf. Wksp. (WCNCW)}, Barcelona, Spain, Apr. 2018.

\bibitem{Wind_data}
\BIBentryALTinterwordspacing
{National Centre for Earth Observation and National Centre for Atmospheric
  Science}, ``{The CEDA Archive: The Natural Environment Research Council's
  Data Repository for Atmospheric Science and Earth Observation}.'' [Online].
  Available: \url{http://archive.ceda.ac.uk/}
\BIBentrySTDinterwordspacing

\bibitem{Turbine_data}
\BIBentryALTinterwordspacing
{The Kingfisher Information Service - Offshore Renewable Cable Awareness},
  ``{Awareness Chart of Walney 1-4}.'' [Online]. Available:
  \url{http://www.kis-orca.eu/downloads}
\BIBentrySTDinterwordspacing

\bibitem{2001-UAV-energy-model}
G.~D. A.~R. S.~Bramwell and D.~Balmford, \emph{{Bramwell’s Helicopter
  Dynamics, 2nd ed.}}\hskip 1em plus 0.5em minus 0.4em\relax American Institute
  of Aeronautics \& Ast (AIAA), 2001.

\bibitem{2006-UAV-energy-model}
A.~Filippone, \emph{{Flight Performance of Fixed and Rotary Wing
  Aircraft}}.\hskip 1em plus 0.5em minus 0.4em\relax American Institute of
  Aeronautics \& Ast (AIAA), 2006.

\end{thebibliography}
\end{footnotesize}}

\end{document}